\documentclass[11pt] {article}
\usepackage{authblk}
\usepackage{etex}
\usepackage[utf8]{inputenc}
\usepackage[english]{babel}
\usepackage{amsmath}
\usepackage{amsthm}
\usepackage{amssymb}
\usepackage{sectsty}
\usepackage{titlesec}
\usepackage{color}
\usepackage[color,matrix,arrow]{xy}
\usepackage{amsgen}
\usepackage{amstext}
\usepackage{amsbsy}
\usepackage{amsopn}
\usepackage{amsfonts}
\usepackage{eepic}
\usepackage{graphicx}
\usepackage{epsf}
\usepackage{mathtools}
\usepackage{pstricks}
\usepackage{float}
\usepackage{enumerate}
\usepackage{enumitem}
\usepackage{tikz}
\usepackage{faktor}
\usepackage{tikz-cd}
\usetikzlibrary{decorations.pathmorphing}
\xyoption{all}
\usetikzlibrary{patterns}

\usepackage{pgf}
\usepackage{tikz}

\newcommand{\footrecall}[1]{}

\titleformat*{\section}{\large\bfseries}
\titleformat*{\subsection}{\normalsize \bfseries}

\newcommand{\N}{\mathbb{N}}
\newcommand{\ZZ}{\mathbb{Z}}
\newcommand{\Z}{\mathbb{Z}}

\newcommand{\RR}{\mathbb{R}}
\newcommand{\R}{\mathbb{R}}

\newcommand{\Fix}{\text{Fix}}

\newcommand{\Ker}{\text{Ker}}
\newcommand{\End}{\text{End}}
\newcommand{\Inn}{\text{Inn}}
\newcommand{\Aut}{\text{Aut}}

\newcommand{\Img}{\text{Im}}
\newcommand{\Haus}{\text{Haus}}
\newcommand{\diam}{\text{diam}}

\newcommand{\mc}{\mathcal}

\theoremstyle{definition}
\newtheorem{Thm}{Theorem}[section]
\newtheorem{Cor}[Thm]{Corollary}

\newtheorem{Prop}[Thm]{Proposition}
\newtheorem{Lem}[Thm]{Lemma}
\newtheorem{Ex}[Thm]{Example}
\newtheorem{Rmk}[Thm]{Remark}

\newtheorem{Question}[Thm]{Question}
\newtheorem*{theorem*}{Theorem}

\title{On endomorphisms of automatic groups}

\author{Andr\'e Carvalho \thanks{andrecruzcarvalho@gmail.com}}
\affil{Centre of Mathematics, University of Porto, R. Campo Alegre, 4169-007 Porto,
Portugal}

\date{}

\begin{document}

% --------------------------------------------------------------
%                         Começar aqui
% --------------------------------------------------------------
 
%\title{Stuff}

\maketitle

\begin{abstract}
We propose two geometric versions of the bounded reduction property and find conditions for them to coincide. In particular, for the natural automatic structure on a hyperbolic group, the two notions are equivalent. We study endomorphisms with $L$-quasiconvex image and prove that those with finite kernel satisfy a synchronous version of the bounded reduction property. Finally, we use these techniques to prove $L$-quasiconvexity of the equalizer of two endomorphisms under certain (strict) conditions.
\end{abstract}

\section{Introduction}
The study of fixed subgroups of endomorphisms of groups started with the (independent) work of Gersten \cite{[Ger87]} and Cooper \cite{[Coo87]}, using respectively graph-theoretic and topological approaches. They proved that the subgroup of fixed points $\Fix(\varphi)$ of some fixed automorphism $\varphi$ of $F_n$ is always finitely generated, and Cooper succeeded on classifying from the dynamical viewpoint the fixed points of the continuous extension of $\varphi$ to the boundary of $F_n$. Bestvina and Handel subsequently developed the theory of train tracks to prove that $\Fix(\varphi)$ has rank at most $n$ in \cite{[BH92]}. The problem of computing a basis for $\Fix(\varphi)$ had a tribulated history and was finally settled by Bogopolski and Maslakova in 2016 in \cite{[BM16]}.
This line of research extended early to wider classes of groups. For instance, Paulin proved in 1989 that the subgroup of fixed points of an automorphism of a hyperbolic group is finitely generated \cite{[Pau89]}. Fixed points were also studied for right-angled Artin groups \cite{[RSS13]} and lamplighter groups \cite{[MS18]}.
 
One of the essential tools used in proving these results is the bounded reduction property (also known as the bounded cancellation lemma) introduced in \cite{[Coo87]} and followed by many others. 
In \cite{[Car21]}, the author provided three equivalent geometric definitions of this property for hyperbolic groups and proved that uniformly continuous endomorphisms with respect to a visual metric are precisely the ones with finite kernel and quasiconvex image.

Automatic groups were introduced in \cite{[ECHLPT92]} and this class of groups is very large and so, it can be hard to deal with. In particular, it contains direct products of hyperbolic groups. Despite being an interesting class with many developments and the fact that study of endomorphisms of groups plays an important role in the theory of finitely generated groups, not much has been done regarding the study of endomorphisms of automatic groups.

In this paper, we will adapt a geometric definition of the bounded reduction property proposed in \cite{[Car21]} to the class of automatic groups and consider another version of this property, finding conditions for both versions to coincide. In particular, if the group is hyperbolic, the two versions are equivalent, extending in this way the results from \cite{[Car21]}. We also propose a synchronous version of the BRP and prove some technical results that help us prove the BRP for some endomorphisms. 

In Section \ref{secqci}, we focus on endomorphisms with $L$-quasiconvex image and generalize a well-known result on hyperbolic groups to the class automatic groups:

\newtheorem*{estsync}{Theorem \ref{estsync}}
\begin{estsync}
Let $G$ be an automatic group with an automatic stucture $L$ for $\pi_1:A^*\to G$ and $\varphi$ be a virtually injective endomorphism with $L$-quasiconvex image. Then there is some automatic structure $K$ such that the synchronous BRP holds for $(\varphi, K, L)$.
\end{estsync}

In Section \ref{secapps} we apply the techniques developed in the previous section to obtain some finiteness results on equalizers (and so, on fixed subgroups) of endomorphisms of automatic groups:
 
\newtheorem*{eq}{Corollary \ref{eq}}
\begin{eq}
Let $G_1$ and $G_2$ be automatic groups with automatic structures $L_1$ and $L_2$ for $\pi_1:A^*\to G_1$ and $\pi_2:B^*\to G_2$, respectively. Let $\varphi,\psi:G_1\to G_2$ be  homomorphisms such that the synchronous BRP holds for $(\varphi, L_1,L_2)$ and $(\psi, L_1,L_2)$. 
Then $\text{Eq}(\varphi,\psi)=\{x\in G_1\mid x\varphi=x\psi\}$ isomorphic to a $(L_2\diamond L_2)$-quasiconvex subgroup of $G_2\times G_2$. In particular, $\text{Eq}(\varphi,\psi)$ is automatic.
\end{eq}

Even though the hypothesis are quite strong, as we will remark in  Section \ref{secapps} proving exactly for which endomorphisms they hold in the case where the group is free and the structure considered is the structure of all geodesics, this result provides an alternative proof of a  result from \cite{[GS91]} which concerns the quasiconvexity of the centralizer of a finite subset:

\newtheorem*{centralizer}{Corollary \ref{centralizer}}
\begin{centralizer}
 The centralizer of a finite subset of a biautomatic group is biautomatic.
\end{centralizer}

We remark that although the results are obtained under strict hypothesis, very little was known about endomorphisms of automatic groups.

The paper is organized as follows. In Section \ref{secpreliminares}, we present some preliminaries on automatic groups. In Section \ref{secbrp}, we present two versions of the bounded reduction property for automatic groups and find some relations between them.
In Section \ref{secqci} we extend some results obtained in the previous section and focus on endomorphisms with an $L$-quasiconvex image. In Section \ref{secapps}, we apply the techniques developed in the paper to prove quasiconvexity of some equalizers of endomorphisms. In Section \ref{secclass}, we show a visual representation of the properties studied in this paper for endomorphisms of (virtually) free, hyperbolic and automatic groups. 
We finish with some questions in Section \ref{secquestions}.

\section{Preliminaries}
\label{secpreliminares}
We now introduce some known results on automatic groups. For more details, the reader is referred to  \cite{[ECHLPT92]} and \cite{[HRR17]}.
Let $G$ be a group, $A$ be a finite alphabet and $\pi:A^*\to G$ be a surjective homomorphism. A language $L\subseteq A^*$ is said to be a \emph{section of $\pi$} if $L\pi=G$ and if $\pi\lvert_L$ is bijective, we say that $L$ is a \emph{transversal of $\pi$}. Let $\$$ be a symbol that doesn't belong to $A$ and consider 
\[A_\$=(A\times A)\cup (A\times \{\$\})\cup (\{\$\}\times A).\]
Given two words $u, v \in A^*$, the convolution $u \diamond v$ is the only word in $A_{\$}^*$ such that the projection on the first (respectively second) component belongs to $u\$^*$ (respectively $v\$^*$).

The language $L$ is an \emph{automatic structure for $\pi$} if it is a rational section of $\pi$ and, for every $x \in A\cup\{1\}$, there are finite state automata $M_x$ over $A_\$$ such that $L(M_x)=\{u\diamond v \mid u,v\in L,\, v\pi=(ux)\pi\}$. If, additionally, $L$ is a transversal, then we say that  $L$ is an \emph{automatic structure with uniqueness for $\pi$}. Moreover, if there are finite state automata ${}_{x}M$ over $A_\$$ such that $L({}_{x}M)=\{u\diamond v \mid u,v\in L,\, v\pi=(xu)\pi\}$, $L$ is a \emph{biautomatic structure for $\pi$}.

It is well known that, if $L$ is an automatic structure for $\pi$, then the set of all shortlex minimal representatives of elements of $G$ in $L$ form an automatic structure with uniqueness for $\pi$. 

 Given a word $u\in \widetilde{A}^*$, we denote by $u^{[n]}$ the prefix of $u$ with $n$ letters. If $n>\lvert u\lvert $, then we consider $u^{[n]}=u$. Given a group $G=\langle A\rangle$, consider its Cayley graph $\Gamma_A(G)$ with respect to $A$ endowed with the \emph{geodesic metric} $d_A$, defined by letting $d_A(x,y)$ to be the length of the shortest path in $\Gamma_A(G)$ connecting $x$ to $y$. We will slightly abuse notation when we write $d_A$. Indeed, when we are given an automatic structure $L$ for $\pi:A^*\to G$, we have that $G=\langle A\pi \rangle$ and we will write $d_A$ to denote $d_{A\pi}$. 
 
 For a language $L\subseteq A^*$ we say that the \emph{fellow traveller property holds for $L$} if there is some $N\in \mathbb N$ such that, for all $u,v\in L$, 
 \[d_A(u\pi,v\pi)\leq 1\Rightarrow d_A(u^{[n]}\pi,v^{[n]}\pi)\leq N,\]
  for every $n\in \N$.  

 We say that $u,v \in A^*$ are \emph{$p$-fellow travellers}, or that $u$ and $v$ \emph{$p$-fellow travel}, in $\Gamma_A(G)$ if $d_A(u^{[n]}\pi,v^{[n]}\pi) \leq p$ for every $n \in \mathbb{N}$.

We will consider similar definitions to the ones above when the paths remain at bounded distance but asynchronously. We say that two words $u$ and $v$ in $L$ \emph{asynchronously $p$-fellow travel} if they $p$-fellow travel up to some reparametrisations, meaning that there are nondecreasing surjective functions $\phi,\psi:\N\to \N$ such that for all $t\in \N$ we have that $d_A\left(u^{[\phi(t)]}\pi,v^{[\psi(t)]}\pi\right)\leq p.$ We say that a rational section $L$ is an \emph{asynchronous automatic structure} if the asynchronous version of the fellow traveller property holds for $L$.

Our definition coincides with the one used in \cite{[HHRS21]}. It is slightly different from the most common definition (see, for example \cite{[HRR17]},\cite{[NS94]},\cite{[Sha92]}) because we take functions from $\N$ to $\N$ instead of taking the paths parametrised by arc length and nondecreasing surjective maps $\phi,\psi: \RR_{\geq 0}\to  \RR_{\geq 0}$ such that for all $t\in \RR_{\geq 0}$ we have that $d_A\left(u(\phi(t)),v(\psi(t))\right)\leq p.$ These two definitions can be seen to be equivalent. Indeed, if $L$ is an asynchronous automatic structure with the discrete definition, by linear interpolation, we can prove that it is also an asynchronous automatic structure for the continuous version of the definition. Conversely, suppose that $L$ is an asynchronous automatic structure for the continuous version of the definition. Let $u,v\in L$ ending at distance at most one apart and take  $\phi,\psi:[0,\infty)\to [0,\infty)$ such that $$d(u(\phi(t)),v(\psi(t)))\leq C,$$
for all $t\in [0,\infty)$. Let $t_i$ be such that $\phi(t_i)=i$. Notice that we can take $t_0=\phi(t_0)=\psi(t_0)=0$.
Define $\phi',\psi':\N\to \N$ as:

For $ \lfloor \psi(t_i) \rfloor+i \leq n < \lfloor \psi(t_{i+1}) \rfloor+i+1$, let $\phi'(n)=i$ and $\psi'(n)=n-i$.

We have that $\phi'$ and $\psi'$ are nondecreasing. Indeed in each interval of the form $ \lfloor \psi(t_i) \rfloor+i \leq n < \lfloor \psi(t_{i+1}) \rfloor+i+1$, $\phi'$ is constant and increases by $1$ when we move to the next interval and $\psi'$ increases by $1$ in each step inside the interval and remains constant when we change interval. Clearly, they are both surjective as well.

 Now let $n\in \N$ and take $i\in \N$ such that $ \lfloor \psi(t_i) \rfloor+i \leq n < \lfloor \psi(t_{i+1}) \rfloor+i+1$. 
 Then, $\lfloor \psi(t_i) \rfloor\leq n-i< \lfloor \psi(t_{i+1}) \rfloor+1$, and since $n-i$ is an integer, then,  
 \begin{align}
 \label{c1}
 \lfloor \psi(t_i) \rfloor\leq n-i \leq \lfloor \psi(t_{i+1}) \rfloor.
 \end{align}
Now, if  $ \psi(t_i) \leq n-i$, there is some $t\in [t_i,t_{i+1}]$ such that $\psi(t)=n-i$. Fix such a $t$. Since $\phi$ is nondecreasing, then $i=\phi(t_i)\leq \phi(t)\leq \phi(t_{i+1})=i+1$. It follows that 
 \begin{align*}
d(u^{[\phi'(n)]},v^{[\psi'(n)]})&=  d(u^{[i]},v^{[n-i]})\\
&= d(u^{[i]},v(\psi(t)))\\
&\leq d(u^{[i]},u(\phi(t)))+d(u(\phi(t)),v(\psi(t)))\\
& \leq  d(u(\phi(t_i)),u(\phi(t)))+C\\
&\leq C+1.
\end{align*}

If $n-i<\psi(t_i) $, then, by (\ref{c1}), $\psi(t_i)-(n-i)<1$. In this case, let $t$ be such that $\psi(t)=n-i$. We have that $\psi(t_i)-\psi(t)<1$. Then 
 \begin{align*}
d(u^{[\phi'(n)]},v^{[\psi'(n)]})&=  d(u^{[i]},v^{[n-i]})\\
&= d(u^{[i]},v(\psi(t)))\\
&\leq d(u(\phi(t_i)),v(\psi(t_i)))+d(v(\psi(t_i)),v(\psi(t)))\\
& \leq  C+1.
\end{align*}

There is also another common way to define an asynchronous automatic structure (see, for example, \cite{[NS92]}, \cite{[NS95]}). According to this definition, $L$ is an asynchronous automatic structure if there is some $p>0$ such that for all words $u$ and $v$ in $L$  ending at distance at most $1$, there is a nondecreasing surjective function $\phi:\RR_{\geq 0}\to \RR_{\geq 0}$ such that for all $t\in \RR_{\geq 0}$, we have that 
\begin{align}
\label{asincshap}
d_A\left(u(t),v\left(\phi(t)\right)\right)\leq p.
\end{align}
This can be easily seen to be equivalent to the previous definitions using \cite[Lemma 2.1]{[Neu92]}. We remark that this notion of an asynchronous automatic structure, when considered in its discrete version does not coincide with the continuous version (nor with the previous definition). Similarly, using linear interpolation, we can see that if some $L$ verifies the discrete condition, then it verifies the continuous one, but the converse is not true. For example, if we take $\ZZ$ generated by $A=\{a,a^{-1},a^2\}$, $\pi:A^*\to \ZZ$ defined in the natural way and put $L=a^*\cup (a^{-1})^*\cup (a^2)^*$, then $L$ is rational and $\pi\lvert_L$ is surjective. Let $u,v\in L$ be such that $d_A(u,v)\leq 1$. In particular, we must have that $u\pi \times v\pi\geq 0$. If $u\pi,v\pi\leq 0$, then it is obvious that $d(u(t),v(t))\leq 1$, for all $t\in[0,\infty)$. If either $u\pi$ or $v\pi=0$, then we also have that $d(u(t),v(t))\leq 1$, for all $t\in[0,\infty)$. So, suppose that  $u\pi,v\pi> 0$. If $u=a^{k_1}$ and $v=a^{k_2}$ with $k_1,k_2\in \N$, then $\lvert k_1-k_2\lvert\, \leq 2$ and $u$ and again  $d(u(t),v(t))\leq 1$, for all $t\in[0,\infty)$. Similarly, if  $u=(a^2)^{k_1}$ and $v=(a^2)^{k_2}$, with $k_1,k_2\in \N$, then $d(u(t),v(t))\leq 1$, for all $t\in[0,\infty)$. If $u=(a^2)^{k_1}$ and $v=a^{k_2}$ with $k_1,k_2\in \N$, then, $\lvert 2k_1-k_2\lvert\,\leq 1$ and letting $\phi(x)=2x$, we have that $d(u(t),v(\phi(t)))\leq 1$. Similarly, if $u=a^{k_1}$ and $v=(a^2)^{k_2}$, we get the same with $\phi(x)=\frac 1 2 x$. So, it is an asynchronous automatic structure  according to (\ref{asincshap}).
However, in a discrete version of (\ref{asincshap}), for every $C\in \N$, if we take $u=(a^2)^{C+1}$ and $v=a^{2C+2}$, for every  surjective, nondecreasing $\phi:\N\to \N$, we have that $d(u(C+1),v(\phi(C+1)))\geq d(a^{2C+2},a^{C+1})=C+1>C$, because $\phi(n)\leq  n$ for all $n\in \N$.

 It is also a well-known fact that a rational section $L$ is an automatic structure for $\pi$ if and only if the fellow traveller property holds for $L$. Similarly, biautomatic structures can be characterized by a variant of the fellow traveller property, where points might start at distance at most one apart.  If $L_1\subseteq A^*$ and $L_2\subseteq B^*$ are synchronous or asynchronous automatic structures, then $L_1\cup L_2$ is rational on $(A\cup B)^*.$ We say that $L_1$ and $L_2$ are \emph{equivalent} if $L_1\cup L_2$ is an asynchronous automatic structure. For more information on equivalent automatic structures the reader is referred to \cite{[NS92]}. We remark that this is indeed an equivalence relation and that transitivity is easier to check using the definition of an asynchronous automatic structure in (\ref{asincshap}) and taking the composition of the reparametrisation mappings.
In hyperbolic groups, all structures are equivalent.  We will also consider the notion of \emph{synchronous equivalence} when $L_1\cup L_2$ is a synchronous automatic structure. This notion is stricter than the notion of equivalence. For example in $\ZZ$, the structure given by $(aa^{-1}a)^*\cup (a^{-1}aa^{-1})^*$ is not synchronous equivalent to the usual structure given by $a^*\cup (a^{-1})^*$. In fact, taking $\ZZ$ with a single generator $a=1$, there are infinitely many synchronous equivalence classes. Indeed, for $n\in \N$, the structures $L_n=(a(a^{-1}a)^n)^*\cup (a^{-1}(aa^{-1})^n)^*$ belong to different equivalence classes.

Let $(X,d)$ be a metric space and $Y,Z$ nonempty subsets of $X$. We call the $\varepsilon$-neighbourhood of $Y$ in $X$ and we denote by $\mc V_\varepsilon(Y)$ the set $\{x\in X\mid d(x,Y)\leq \varepsilon\}$.
The {\em Hausdorff distance} between
$Y$ and $Z$, which we denote by Haus$(Y,Z)$, is given by
$$\Haus(Y,Z)=\inf\{ \varepsilon>0\mid Y\subseteq \mc V_\varepsilon (Z) \text{ and } Z\subseteq\mc V_\varepsilon(Y)\}$$
if the set is nonempty. Otherwise Haus$(Y,Z)=\infty$.

We say that an automatic structure has the \emph{Hausdorff closeness property}  if there is some $N\in \mathbb N$ such that, for all $u,v\in L$ such that $d_A(u\pi,v\pi)\leq 1$, the paths in the Cayley graph defined by $u$ and $v$ are at Hausdorff distance bounded by $N$.
 
 We will also use a characterization of a \emph{boundedly asynchronous automatic structure} (Theorem 7.2.8 in \cite{[ECHLPT92]}). We will only use the result to prove that a certain structure is asynchronous automatic and may take this characterization as a definition. 
 
 If $G$ is an automatic group with automatic structure $L$ for $\pi:A^*\to G$, a \emph{departure function} for $(G,L)$ is any function $D:\R_0^+\to \R_0^+$ such that, if $w\in L$, $r,s\geq 0$, $t\geq D(r)$ and $s+t\leq |w|$, then $d_A(w^{[s]}\pi,w^{[s+t]}\pi)>r.$
 \begin{Lem}
 \label{async}
 Let $G$ be a group, $A$ a set of semigroup generators for $G$ and $L$ a regular language over $A$ that maps onto $G$. Then $L$ is a boundedly asynchronous structure on $G$ if and only if the following conditions are satisfied.
 \begin{enumerate}
 \item there exists a departure function $D$ for $(G,L)$; and
 \item there exists a constant $K>0$ such that, for every pair of strings $w,w'\in L$ whose images under the map $\pi:A^*\to G$ are a distance at most one apart in $\Gamma_A(G)$, the paths defined by $w$ and $w'$ are at most Hausdorff distance $K$ from each other.
 \end{enumerate}
 \end{Lem}
 
 While an automatic structure is obviously asynchronous automatic and satisfies condition 2 of the above Lemma, it is not true that it always admits a departure function. For example, take a finite group $G$ and the alphabet $A=G$. Then take $L=A^*$ and $\pi:A^*\to G$. Let $p=\diam(G)$. It is obvious that any two words $p$-fellow travel, thus $L$ is an automatic structure. However, it is easy to see that it can never admit a departure function. If there was one, say $D:\R_0^+\to \R_0^+$, then we could take $a\in A$ and the word $w=a^{D(0)|A|}$ and taking the prefixes of size a multiple of $D(0)$, $a^{kD(0)}$, for $k=0,\ldots, A$, then we must have two distinct ones, say $a^{iD(0)}$ and $a^{jD(0)}$, with $0\leq i<j\leq |A|$ representing the same element in $G$. This means that for $r=0$, $s=iD(0)$, $t=(j-i)D(0)>D(0)$, we have that $s+t\leq|w|$ and $d_A(w^{[s]}\pi,w^{[s+t]}\pi)=0$.
 
 However, if $L$ is an automatic structure, there is some $L'\subseteq L$ such that $L'$ is boundedly asynchronous automatic. Also, every boundedly asynchronous automatic structure is in particular an asynchronous automatic structure.

We say that an automatic structure is \emph{factorial} if and only if  it is closed under taking factors, i.e., given $u\in L$, if $v$ is a factor of $u$, then $v\in L$. 

If a group $G$ admits a factorial (bi)automatic structure, then we say that $G$ is \emph{f-(bi)automatic}.  We can see that this is not a very strict condition to impose, following closely the proof of Theorem 2.5.9 in \cite{[ECHLPT92]}.

\begin{Prop}
Let $G$ be a biautomatic group with a biautomatic structure  $L$ for $\pi:A^*\to G$. Then Pref($L$) is also a biautomatic structure of $G$.
\end{Prop}
\noindent\textit{Proof.}  It is well know that Pref($L$) is rational and since $L\subseteq$ Pref($L$), then $\pi|_{\text{Pref}(L)}$ is surjective. We now verify that the biautomatic variant of the fellow traveller property holds for Pref($L$). Let $u,v\in$ Pref($L$) be such that there are $x,y\in A\cup\{\varepsilon\}$ for which $v\pi=(xuy)\pi$. Let $c$ be the number of states of a finite state automaton accepting $L$. Let $u'$ and $v'$ be words of length at most $c$ such that $uu',vv'\in L$. We have that $d_A((xuu')\pi,(vv')\pi)\leq d_A((xuu')\pi,(xu)\pi)+d_A((xu)\pi,v\pi)+d_A(v\pi,(vv')\pi)\leq 2c+1.$ Since $L$ is biautomatic, we have that $xuu'$ and $vv'$ $N$-fellow travel where $N$ is a constant depending only on $c$. Thus, $xu$ and $v$ $(N+2c)$-fellow travel. 
\qed

\begin{Prop}
\label{bifatorial}
Every biautomatic group is f-biautomatic.
\end{Prop}
\noindent\textit{Proof.} Let $G$ be a biautomatic group with a biautomatic structure $L$ for $\pi:A^*\to G$ and let $c$ be the number of states of a finite automaton recognizing $L$ and $N$ be the biautomatic fellow traveller constant for $L$. Then, Pref($L$) is a biautomatic structure. Now, we claim that Suff$(L)$ is also a biautomatic structure. Indeed, let $u,v\in \text{Suff($L$)}$ such that there is $a\in A\cup \{\varepsilon\}$ for which $d_A((au)\pi=v\pi)\leq 1$.
  Take words $u',v'\in A^*$ of length less than $c$ such that $u'u$ and $v'v$ belong to $L$.  Then $|v'a(u')^{-1}|\leq 2c+1$ and 
  
  \[d_A((v'a(u')^{-1}u'u)\pi,(v'v)\pi)=d_A(au ,v)\leq 1\] and so $u'u$ and $v'v$ are paths represented by words in $L$ that start at bounded distance and end at distance $1$. By the biautomatic fellow traveller property of $L$, they $k(2c+1)$-fellow travel  and so $u$ and $v$ $(2c+k(2c+1))$-fellow travel. 
  So, Pref(Suff($L$)) is a factorial biautomatic structure.
  \qed\\

In particular, every hyperbolic group admits a factorial biautomatic structure with uniqueness. Indeed,  the set of geodesics $Geo_A(G)$ forms  such a structure. 
If we take the natural automatic structure with uniqueness given by the set of  all shortlex minimal representatives of elements of $G$ in $Geo_A(G)$, we have a biautomatic structure with uniqueness. It must also be factorial, because for $u\in Geo_A(G)$, letting $\overline u$ denote the shortlex minimal word such that $\overline u \pi=u\pi$ and taking a factor $v$ of $\overline u$, then $\overline u=wvw'$ for some words $w,w'\in A^*$. If $v\neq \overline v$, then $w\overline v w'\in Geo_A(G)$ is such that $w\overline v w'<_{sl} \overline u$ and $(w\overline vw')\pi=u\pi$, which is absurd.

Given two automatic groups $G_1$ and $G_2$ with automatic structures with uniqueness $L_1$ and $L_2$, respectively, then $L_1L_2$ is an automatic structure with uniqueness for $G_1\times G_2$, so, the product of two hyperbolic groups is f-biautomatic and admits a factorial biautomatic structure with uniqueness.
Moreover, a factorial biautomatic structure with uniqueness must admit a departure function. Let $D:\R_0^+\to \R_0^+$ be  defined by $x\mapsto \max\{|w|  \;\big{|}\; w \in L,  d_A(1,w\pi)\leq x\}$. Then, given a word $w\in L$,  $r,s\geq 0$, $t\geq D(r)$ and $s+t\leq |w|$, then $d_A(w^{[s]}\pi,w^{[s+t]}\pi)>r$, because the path from $w^{[s]}$ to $w^{[s+t]}$ must be a word in $L$ of length greater than $D(r)$ which means by definition of $D$ that it represents an element with norm greater than $r$.

So, we know that a biautomatic group admits a factorial biautomatic structure and it admits a structure with a departure function. We now prove that we always have a structure with both these properties.

\begin{Prop}
Every biautomatic group has a factorial automatic structure admitting a departure function.
\end{Prop}

\noindent\textit{Proof.}
Let $L\subseteq A^*$ be a biautomatic structure for $\pi:A^*\to G$. Then, it admits a substructure $L'$ with a departure function $D$. Since $L$ is biautomatic, then $L'$ must also be biautomatic. Then Fact($L'$) is biautomatic by Proposition \ref{bifatorial}. We claim that $D$ is also a departure function for Fact($L'$). Let $w\in \text{Fact($L'$)}$ and $r,s\geq 0$ and $t\geq D(r)$ such that $s+t\leq |w|$. Then there are $u,v\in A^*$ such that $uwv\in L'$. We have that 
\[d_A(w^{[s]}\pi,w^{[s+t]}\pi)=d_A((uw^{[s]})\pi,(uw^{[s+t]})\pi)=d_A((uwv)^{[|u|+s]}\pi,(uwv)^{[|u|+s+t]}\pi)>r.\]
\qed

When we work with an automatic structure, given an element $g\in G$, we will denote by $\overline g$ an arbitrary minimal length element  of $L$ such that $\overline g\pi=g$.

The following lemma, known as the bounded length difference lemma, will be crucial to our work.
\begin{Lem}[\cite{[ECHLPT92]}, Lemma 2.3.9]
\label{kuniq}
Let $G$ be an automatic group, $A$ be an alphabet, $\pi:A^*\to G$ be a surjective homomorphism and $L$ be an automatic structure $\pi$. Then, there is a constant $K>0$ such that, if $w\in L$ and $g\in G$ is a vertex of the Cayley graph at distance at most one from $w\pi$, we have the following situation.
\begin{enumerate}
\item $g$ has some representative of length at most $|w|+K$ in L; and 
\item if some representative of $g$ in $L$ has length greater that $|w|+K$, there are infinitely many representatives of $g$ in $L$.
\end{enumerate} 
\end{Lem}

As a corollary, we have that if $g$ is at bounded distance from $w\pi$, then there is some representative of bounded length.
In particular, it follows that, for all $g,h\in G$, \[\left | |\bar g|-|\bar h |\right|\leq Kd_A(g,h).\]

If $L$ is an automatic structure for $\pi:A^*\to G$ such that the empty word $\varepsilon$ does not belong to $L$, we can consider $L'=L\cup \{\varepsilon\}$, which is obviously a rational section and satisfies the fellow traveler property. Thus,  for all $g\in G$, we have that 
\begin{align}
\label{minrep}
|\bar g|\leq Kd_A(1,g).
\end{align}

 A subgroup $H\leq G$ is said to be $L$-quasiconvex if there exists some $N$ such that, for all $h,h'\in H$ and $w\in L$ such that $h(w\pi)=h'$, we have that $$d_A(h(w^{[n]}\pi),H)\leq N,$$
holds for all $n\in \N$.

\section{Bounded Reduction Property}
\label{secbrp}
The purpose of this section is to explore the notion of bounded reduction for endomorphisms of automatic groups. We will 
adapt a definition proposed in \cite{[Car21]} and propose a new one, establishing some relations between them.

We start by presenting some definitions concerning hyperbolic groups.

 A {\em quasi-isometric embedding} of metric spaces is a mapping $\varphi:(X,d) \to
(X',d')$ such that there exist constants $\lambda \geq 1$ and $K \geq 0$
satisfying
$$\frac{1}{\lambda}d(x,y) -K \leq d'(x\varphi,y\varphi) \leq \lambda d(x,y) +K$$
for all $x,y \in X$. We may call it a
$(\lambda,K)$-quasi-isometric embedding if we want to stress the
constants. 
 A \emph{$(\lambda,K)$-quasi-geodesic} of between $x$ and $y$ is a $(\lambda,K)$-quasi-isometric embedding $\xi : [0, s] \to X$ such that $0\xi = x$
and $s\xi = y$, where $[0, s] \subset \RR$ is endowed with the usual metric of $\RR$.

In a hyperbolic group $H$, given $g,h,p\in H$, we define the \emph{Gromov product of $g$ and $h$} taking $p$ as basepoint by
$$(g\lvert h)_p^A=\frac{1}{2} (d_A(p,g)+d_A(p,h)-d_A(g,h)).$$
We will write $(g\lvert h)$ to denote $(g\lvert h)_1^A$.

In \cite{[Car21]}, several geometric versions of the Bounded Reduction Property are introduced. We say that the \emph{hyperbolic BRP holds} for a hyperbolic group endomorphism if one of the (equivalent) conditions in the following theorem hold.
\begin{Thm}
Let $H$ be a hyperbolic group and $\varphi\in \End(H)$. The following conditions are equivalent:
\label{equivsbrphiper}
\begin{enumerate}[label=\roman*.]
\item for every $p\geq 0$ there is some $q\geq 0$ such that:
given two geodesics $u$ and $v$ such that  the concatenation $u+v$ is a $(1,p)$-quasi-geodesic,  we have that given any two geodesics $\alpha$, $\beta$, from $1$ to $u\varphi$ and from $u\varphi$ to $(uv)\varphi$, respectively, the concatenation $\alpha+\beta$
is a $(1,q)$-quasi-geodesic.

\item there is some $q\geq 0$ such that:
given two geodesics $u$ and $v$ such that the concatenation $u+v$ is a geodesic,  we have that given any two geodesics $\alpha$, $\beta$, from $1$ to $u\varphi$ and from $u\varphi$ to $(uv)\varphi$, respectively, their concatenation $\alpha+\beta$ is a $(1,q)$-quasi-geodesic.
\item $\forall\, p>0\,\exists\, q>0\,\forall\, u,v\in H\, ((u|v)\leq p\Rightarrow (u\varphi|v\varphi)\leq q)$.
\item $\exists\, q>0\,\forall\, u,v\in H\, ((u|v)=0\Rightarrow (u\varphi|v\varphi)\leq q)$.
\item there is some $N\in \N$ such that, for all $x,y\in H$ and every geodesic $\alpha=[x,y]$, we have that $\alpha\varphi$ is at bounded Hausdorff distance to every geodesic $[x\varphi,y\varphi]$.
\item there is some $N\in \N$ such that, for all $x,y\in H$ and every geodesic $\alpha=[x,y]$, we have that $\alpha\varphi\subseteq \mc V_N(\xi)$ for every geodesic $\xi=[x\varphi,y\varphi]$.
\item $\varphi$ is coarse-median preserving.
\end{enumerate}
\end{Thm}

Let $G$ be an automatic group, $A$ be a finite alphabet and $\pi:A^*\to G$ be a surjective homomorphism. Let $L\subseteq A^*$ be an automatic structure for $\pi$ and $\varphi\in \End(G)$.  Whenever we are able to, we will try to state the results as generally as possible, avoiding restrictions to the class of automatic structures that the group admits. 

Given $u\in L$, $x\in G$, let $\theta_u^x:\{0,\ldots, \lvert u\lvert \}\to G$ be the function that maps $n$ to $x(u^{[n]}\pi)$.
 If $A'$ is another set of generators for $G$, we define $$N_{A,A'}=\max (\{d_{A'}(1,a)\mid a\in A\}\cup \{d_A(1,a')\mid a'\in A'\})$$ and it follows that 
 \begin{align}
 \label{difgen}
 \frac 1 {N_{A,A'}}d_{A'}(g,h)\leq d_{A}(g,h)\leq N_{A,A'}d_{A'}(g,h)
 \end{align}
holds for all $g,h\in G.$

Given a homomorphism $\varphi:G_1\to G_2$ between two automatic groups $G_1$ and $G_2$ with automatic structures $L_1\subseteq A^*$ and $L_2\subseteq B^*$ for $\pi_1:A^*\to G$ and $\pi_2:B^*\to G$, respectively, we say that the \emph{BRP holds for $(\varphi,L_1,L_2)$} if there is some $N>0$ such that
$$\Haus\left(\Img(\theta_\alpha^x\varphi),\Img(\theta_{\beta}^{x\varphi})\right)\leq N,$$
for all $x\in G$, $\alpha\in L_1$ and $\beta\in L_2$  such that $\beta\pi_2=\alpha\pi_1\varphi.$ Notice that this Hausdorff distance is not well defined since it is not clear what metric it refers to. However, that distinction is not relevant because of (\ref{difgen}). We will also consider a synchronous version of the BRP. We say that the \emph{synchronous BRP holds for $(\varphi,L_1,L_2)$} if there is some $N>0$ such that
$$d\left((\theta_\alpha^x(n))\varphi,\theta_{\beta}^{x\varphi}(n)\right)\leq N,$$ 
for all $n\in \N$, $x\in G$, $\alpha\in L_1$ and $\beta\in L_2$  such that $\beta\pi_2=\alpha\pi_1\varphi.$

We also propose another definition and we will obtain conditions for them to coincide.
We say that  the \emph{fellow traveller bounded reduction property (FT-BRP) holds for $(\varphi,L_1,L_2)$} if  
for every $p\geq 0$ there is some $q\geq 0$ satisfying the following condition:

Let $u_i \in L_1$ and $u'_i \in L_2 \cap u_i\pi_1\varphi\pi_2^{-1}$ for $i = 1,2,3$. If $u_1u_2,u_3$ are $p$-MFT in $\Gamma_A(G)$, then $u'_1u'_2,u'_3$ are $q$-MFT in $\Gamma_B(G)$.

We now present a technical lemma.

\begin{Lem}
\label{concatftp}
Let $G$ be an automatic group with an automatic structure $L$ for $\pi:A^*\to G$ admitting a departure function $D$. Given $p_1,p_2>0$, there is $N\in \N$ such that given words $u_1,u_2,v_1,v_2\in L$ with the property that $u_1$ and $v_1$ are $p_1$-MFT and $u_2$ and $v_2$ $p_2$-fellow travel, we have that $u_1u_2$ and $v_1v_2$ $N$-fellow travel. If, additionally, $L$ is biautomatic and $d_A((u_1u_2)\pi,(v_1v_2)\pi)$ is bounded, the claim holds for any $p_1$ (non meeting) fellow traveler words $u_1$ and $v_1$.    
\end{Lem}
\noindent\textit{Proof.} Suppose w.l.o.g that $|u_1|\leq|v_1|$. Then, obviously, for $j\leq|u_1|$, we have that $$d_A((u_1u_2)^{[j]},(v_1v_2)^{[j]})=d_A(u_1^{[j]}\pi,v_1^{[j]}\pi)\leq p_1.$$

Also, since $D$ is a departure function we have that 
\begin{align}
\label{dep1}
|v_1|-|u_1|\leq D(2p_1).
\end{align}

Indeed, suppose that $|v_1|-|u_1|>D(2p_1)$. Then, we have that 
\begin{align*}
d_A\left(v_1^{|u_1|}\pi,v_1^{|v_1|}\pi\right)&\leq d_A\left(v_1^{|u_1|}\pi,u_1^{|v_1|}\pi\right)+ d_A\left(v_1^{|v_1|}\pi,u_1^{|v_1|}\pi\right) \\
&\leq d_A\left(v_1^{|u_1|}\pi,u_1^{|u_1|}\pi\right)+ d_A\left(v_1^{|v_1|}\pi,u_1^{|v_1|}\pi\right) \leq 2p_1,
\end{align*}
and that is absurd by the definition of a departure function.

So, we have (\ref{dep1}). For $|u_1|<j\leq |v_1|$, we have that
\begin{align*}
d_A((u_1u_2)^{[j]}\pi,(v_1v_2)^{[j]}\pi)&\leq d_A((u_1u_2)^{[j]}\pi,u_1\pi)+ d_A(u_1\pi,v_1\pi)+ d_A(v_1\pi,(v_1v_2)^{[j]}\pi)\\
&\leq j-|u_1|+ p_1 + |v_1|-j\\
& \leq - |u_1| + p_1 +|v_1|\\
&\leq p_1 +D(2p_1);
\end{align*}

and for $j>|v_1|$, we have that 
\begin{align*}
d_A\left((v_1v_2)^{[j]}\pi, (u_1u_2)^{[j]}\pi\right)\leq \, &d_A\left((u_1u_2)^{[j]}\pi,(v_1v_2)^{[|v_1|+j-|u_1|]}\pi\right) \\&+d_A\left((v_1v_2)^{[|v_1|+j-|u_1|]}\pi, (v_1v_2)^{[j]}\pi\right)\\
\leq\, &   d_A\left(u_2^{[j-|u_1|]}\pi,v_2^{[j-|u_1|]}\pi\right)+|v_1|-|u_1|\\
\leq \, &p_2+D(2p_1).
\end{align*}
Hence,  $u_1u_2$ and $v_1v_2$ are $\max\{2p_1,p_1+D(2p_1),p_2+D(2p_1)\}$-fellow travelers in $\Gamma_A(G).$

 We remark that the fact that $u_1$ and $v_1$ are meeting was only used to prove that case where $j>|v_1|$. If additionally we have that $L$ is biautomatic and that $d_A((u_1u_2)\pi,(v_1v_2)\pi)$ is bounded, then the same argument works since $d_A(u_1\pi,v_1\pi)$ is bounded by $p_1$ and so the biautomatic variant of the fellow traveller property holds. 
 \qed

Now we prove a useful result that will help us establish the  FT-BRP for some endomorphisms.
\begin{Prop}
\label{umetodos}
Let $G$ be an automatic group admitting an automatic structure $L$ for $\pi:A^*\to G$ and a departure function $D$, and take $p>0$. If $u_1,u_2,u_3\in L$ are such that   $u_1u_2,u_3$ are $p$-MFT,  then there is some $q>0$ such that for every $u_1',u_2',u_3'\in L$, for $i=1,2,3$, $u_i'\pi=u_i\pi$, the inequality 
$$d_A((u_1'u_2')^{[n]}\pi, u_3'^{[n]}\pi)\leq q,$$
holds for all $n\in \N.$
\end{Prop}
\noindent\textit{Proof.} Let $p>0$ and  $u_1,u_2,u_3\in L$ such that   $u_1u_2,u_3$ are $p$-MFT.
Take $u_1',u_2',u_3'\in L$ such that, for $i=1,2,3$, $u_i'\pi=u_i\pi$. There is some constant $N$ such that  $u_i$ and $u_i'$ are $N$-fellow travelers because $L$ is automatic.
By Lemma \ref{concatftp}, we have that $u_1'u_2'$ fellow travels $u_1u_2$ and so, it fellow travels $u_3$, which fellow travels $u_3'$. 
\qed

Now, we prove that these versions of the bounded reduction property are preserved when we compose endomorphisms.

\begin{Lem}
\label{compos}
Let $G_1, G_2$ and $G_3$ be automatic groups with automatic structures $L_1\subseteq A_1^*$, $L_2\subseteq A_2^*$, and $L_3\subseteq A_3^*$ for $\pi_1,\pi_2$ and $\pi_3$, respectively, and $\varphi_1:G_1\to G_2$ and $\varphi_2:G_2\to G_3$ be homomorphisms such that the (FT, synchronous) BRP holds for $(\varphi_1,L_1,L_2)$ and the (FT, synchronous) BRP holds for $(\varphi_2,L_2,L_3)$. Then the (FT, synchronous) BRP holds for $(\varphi_1\varphi_2,L_1,L_3)$.
\end{Lem}

\noindent\textit{Proof.} First we deal with the FT-BRP. Let $r>0$. Then there is some $q_1>0$ given by the  FT-BRP for $(\varphi_1,L_1,L_2)$, when we take $p=r$ and some $q_2>0$ given by the FT-BRP for  $(\varphi_2,L_2,L_3)$ when we take $p=q_1$.
  
So, taking words $u,v\in L_1$ such that 
 for every $w\in L_1$ such that $w\pi_1=(uv)\pi_1$, we have that 
$$d_A((uv)^{[n]}\pi_1, w^{[n]}\pi_1)\leq r$$ 
holds for all $n\in \N$, then for all words $u',v',w'\in L_2$, such that $u'\pi_2=u\pi_1\varphi_1$, $v'\pi_2=v\pi_1\varphi_1$ and $w'\pi_2=(uv)\pi_1\varphi_1$, we have that, 
$$d_A((u'v')^{[n]}\pi_2, w'^{[n]}\pi_2) \leq q_1$$
holds for all $n\in \N$, by application of the $(\varphi_1,L_1,L_2)$-FT-BRP. 

So, take words $u',v'\in L_2$ such that $u'\pi_2=u\pi_1\varphi$, $v'\pi_2=v\pi_1\varphi$ and we have that for every $w'\in L_2$ such that $w'\pi_2=(uv)\pi_1\varphi$, the inequality 
$$d_A((u'v')^{[n]}\pi_2, w'^{[n]}\pi_2) \leq q_1$$
holds for all $n\in \N$ and so, by application of the $(\varphi_2,L_2,L_3)$-FT-BRP, it follows that for all words $u'',v'',w''\in L_3$, such that $u''\pi_3=u'\pi_2\varphi_2=u\pi_1\varphi_1\varphi_2$, $v''\pi_3=v'\pi_2\varphi_2=v\pi_1\varphi_1\varphi_2$ and $w''\pi_3=(u'v')\pi_2\varphi_2=(uv)\pi_1\varphi_1\varphi_2$, we have that, 
$$d_A((u''v'')^{[n]}\pi_3, w''^{[n]}\pi_3) \leq q_2$$
holds for all $n\in \N.$

Now we prove the same for the BRP.   Let  $x\in G$, $\alpha\in L_1$ and $\gamma\in L_3$  such that $\gamma\pi_3=\alpha\pi_1\varphi_1\varphi_2.$ Let $\beta$ such that $\beta\pi_2=\alpha\pi_1\varphi_1$ and put $V_{\varphi_2}=\max\{d_{A_3}(1,a\varphi_2)\mid a\in A_2\}$. Since the BRP holds for $(\varphi_1,L_1,L_2)$, there is some $N_1>0$  for which 
$$\Haus\left(\Img(\theta_\alpha^x\varphi_1),\Img(\theta_{\beta}^{x\varphi_1})\right)\leq N_1.$$
This yields that $$\Haus\left(\Img(\theta_\alpha^x\varphi_1\varphi_2),\Img(\theta_{\beta}^{x\varphi_1}\varphi_2)\right)\leq N_1V_{\varphi_2}.$$
  Since  $\gamma\pi_3=\alpha\pi_1\varphi_1\varphi_2=\beta\pi_2\varphi_2,$ by application of the BRP for $(\varphi_2,L_2,L_3)$, there is some $N_2>0$  for which 
$$\Haus\left(\Img(\theta_\beta^{x\varphi_1}\varphi_2),\Img(\theta_{\gamma}^{x\varphi_1\varphi_2})\right)\leq N_2.$$ Thus, 
$$\Haus\left(\Img(\theta_\alpha^x\varphi_1\varphi_2),\Img(\theta_{\gamma}^{x\varphi_1\varphi_2})\right)\leq N_1V_{\varphi_2}+N_2$$
and the BRP holds for $(\varphi_1\varphi_2,L_1,L_3)$.

Finally, we deal with the synchronous BRP. Let  $x\in G$, $\alpha\in L_1$ and $\gamma\in L_3$ be such that $\gamma\pi_3=\alpha\pi_1\varphi_1\varphi_2.$ Let $\beta$ be such that $\beta\pi_2=\alpha\pi_1\varphi_1$, $n\in \N$ and $N_i$ be the constant given by the synchronous BRP for $(\varphi_i,L_i,L_{i+1})$ for $i=1,2$.

We want to prove that $d_{A_3}(\alpha^{[n]}\pi_1\varphi_1\varphi_2,\gamma^{[n]}\pi_3)$ is bounded. Since the synchronous BRP holds for $(\varphi_1,L_1,L_2)$, we have that  $$d_{A_2}(\alpha^{[n]}\pi_1\varphi_1,\beta^{[n]}\pi_2)\leq N_1$$ and so $$d_{A_3} (\alpha^{[n]}\pi_1\varphi_1\varphi_2,\beta^{[n]}\pi_2\varphi_2)\leq N_1V_{\varphi_2}.$$
By application of the synchronous BRP for $(\varphi_2,L_2,L_3)$, we obtain that $$d_{A_3}(\beta^{[n]}\pi_2\varphi_2,\gamma^{[n]}\pi_3)\leq N_2.$$
Hence, we have that $$d_{A_3}(\alpha^{[n]}\pi_1\varphi_1\varphi_2,\gamma^{[n]}\pi_3)\leq N_1V_{\varphi_2}+N_2.$$
\qed

Now, we see that, for an endomorphism and a factorial automatic structure,  the FT-BRP implies the BRP and find sufficient conditions for the converse to hold.

\begin{Prop}
\label{brphaus}
Let $G$ be an  f-automatic group with a factorial automatic structure $L_1$ and an automatic structure $L_2$.  
Let $\varphi$ be an endomorphism such that the FT-BRP holds for $(\varphi,L_1,L_2)$. Then, the BRP holds for $(\varphi,L_1,L_2).$
\end{Prop}
\noindent\textit{Proof.} Let  $\pi_1:A^*\to G$ and $\pi_2:B^*\to G$ be the surjective homomorphisms for which $L_1$ and $L_2$ are automatic structures, respectively. Since $L_1$ is factorial, then the empty word $\varepsilon \in L_1$. We shall suppose that the empty word $\varepsilon \in B^*$ is also in $L_2$. That is not a restriction, since, if $\varepsilon \not\in L_2$, considering $L_2'=L_2\cup\{\varepsilon\}$, we know that $L_2'$ is an automatic structure for the same surjective homomorphism that $L_2$ and the (FT) BRP holds for $(\varphi,L_1,L_2)$ if and only if it holds for $(\varphi,L_1,L_2')$. 

 Take $r_i$ to be the fellow traveller constants for $L_i$ ($i=1,2$), $M$ to be the constant given by the FT-BRP when $p=r_1$,  $B_\varphi=\max\{d_B(1,a\varphi)\mid a\in A\}$, $K$ given by Lemma \ref{kuniq} for $L_2$. Let $x\in G$, $\alpha\in L_1$ and put $y=x(\alpha\pi_1)$. Take $\xi=\overline{(x^{-1}y)\varphi}\in L_2$ and let $k\in \{1,\ldots, |\xi|\}$.

We start by proving that $d_A(\Img(\theta_\alpha^x\varphi), x\varphi(\xi^{[k]}\pi_2))<KB_\varphi+M$, and so, $\Img(\theta_{\overline{\alpha\pi_1\varphi}}^{x\varphi})\subseteq \mc V_{KB_\varphi+M}(\Img(\theta_\alpha^x\varphi) )$.

Since, in $L_2$, $|\overline{(x\varphi)^{-1}x\varphi}|=|\varepsilon|=0$ and $|\overline{(x\varphi)^{-1}y\varphi}|=|\xi|\geq k$, there is some $n_k\in\{1,\ldots |\alpha|\}$ such that $|\overline{(x\varphi)^{-1}(x(\alpha^{[n_k]}\pi_1))\varphi}|< k $ and  $|\overline{(x\varphi)^{-1}(x(\alpha^{[n_k+1]}\pi_1))\varphi}|\geq k$ in $L_2$. Consider the factorization of $\alpha$ given by 
$$\begin{tikzcd}
x \ar[rr,"\alpha^{[n_k]}"] && x_k \ar[rr,"\alpha_{n_k}"] && y.
\end{tikzcd}$$

Since $L_1$ is factorial, then, $\alpha^{[n_k]},\alpha_{n_k}\in L$. Given $w\in L_1$ such that $w\pi_1=\alpha\pi_1$, we know that $w$ and $\alpha$ are $r_1$-fellow travellers in $\Gamma_{A}(G)$ and since the FT-BRP holds for $(\varphi,L_1,L_2)$, then for the $L_2$-geodesics $\beta=\overline{(x^{-1}x_k)\varphi}$ and $\gamma=\overline{(x_k^{-1}y)\varphi}$, we have that, 
$$d_B((\beta\gamma)^{[n]}\pi_2, \xi^{[n]}\pi_2) \leq M$$
holds for all $n\in \N.$
In particular, $d_B((\beta\gamma)^{[k]}\pi_2,\xi^{[k]}\pi_2)<M$ and so, 
\begin{align}
\label{ineq3}
d_B(x\varphi((\beta\gamma)^{[k]}\pi_2),x\varphi(\xi^{[k]}\pi_2))<M.
\end{align}

Now, put $x_{k+1}=x(\alpha^{[n_k+1]})\pi_1$. We have that $$d_B((x^{-1}x_k)\varphi,(x^{-1}x_{k+1})\varphi)=d_B(x_k\varphi,x_{k+1}\varphi)<B_\varphi$$ and so by Lemma \ref{kuniq}, we have that, in $L_2$, $$\left| |\overline{(x^{-1}x_{k+1})\varphi}|-|\overline{(x^{-1}x_{k})\varphi}|\right|<KB_\varphi.$$ Since, by definition of $x_{k+1}$
, $|\overline{(x^{-1}x_{k+1})\varphi}|\geq k$, then 
\begin{align}
\label{ineq4}
d_B(x_{k}\varphi,x\varphi((\beta\gamma)^{[k]}\pi_2))\leq k-|\overline{(x^{-1}x_{k})\varphi}|\leq |\overline{(x^{-1}x_{k+1})\varphi}|-|\overline{(x^{-1}x_{k})\varphi}|<KB_\varphi.
\end{align}

Combining (\ref{ineq3}) and (\ref{ineq4}), we have that $d_B(x_{k}\varphi,x\varphi(\xi^{[k]}\pi_2))<M+KB_\varphi$.
So, we have that $$\Img(\theta_{\overline{\alpha\pi_1\varphi}}^{x\varphi})\subseteq \mc V_{KB_\varphi+M}(\Img(\theta_\alpha^x\varphi) ).$$

Since every word $\beta\in L_2$ such that $\beta\pi_2=\alpha\pi_1\varphi$ $r_2$-fellow travels $\overline{\alpha\pi\varphi}$ in $\Gamma_B(G)$, we have that  $$\Img(\theta_{\beta}^{x\varphi})\subseteq \mc V_{KB_\varphi+M+r_2}(\Img(\theta_\alpha^x\varphi) ).$$

Now, we will prove that $$\Img(\theta_\alpha^x\varphi)\subseteq \mc V_{M}(\Img(\theta_{\overline{\alpha\pi_1\varphi}}^{x\varphi} ))\subseteq \mc V_{KB_\varphi+M+r_2}(\Img(\theta_{\overline{\alpha\pi_1\varphi}}^{x\varphi} )) .$$
Let $v\in \Img(\theta_\alpha^x)$ and consider the factorization
$$\begin{tikzcd}
x \ar[rr,"\alpha_1"] && v \ar[rr,"\alpha_2"] && y.
\end{tikzcd}$$
Since $L_1$ is factorial, then $\alpha_1,\alpha_2\in L_1.$ Setting $\alpha_1'=\overline{(x\varphi)^{-1}v\varphi}\in L_2$, $\alpha_2'=\overline{(v\varphi)^{-1}y\varphi}\in L_2$ and $\xi'=\overline{(x\varphi)^{-1}y\varphi}\in L_2$, and using the FT-BRP, we have that the concatenation $\alpha_1'\alpha_2'$ and $\xi'$ $M$-fellow travel in $\Gamma_B(G)$. In particular $$d_B(v\varphi,\Img(\theta_{\overline{\alpha\pi_1\varphi}}^{x\varphi} ))<M.$$

Again, since every word $\beta\in L_2$ such that $\beta\pi_2=\alpha\pi_1\varphi$ $r_2$-fellow travels $\xi'$, we have that  $$ \Img(\theta_\alpha^x\varphi) \subseteq \mc V_{M+r_2}(\Img(\theta_{\beta}^{x\varphi})).$$
\qed

\begin{Cor}
\label{strongqc}
Let $G$ be an  f-automatic group with  factorial automatic structure $L_1$ and an automatic structure $L_2$. Let $\varphi$ be an endomorphism such that the FT-BRP holds for $(\varphi,L_1,L_2)$. Then, $G\varphi$ is $L_2$-quasiconvex.
\end{Cor}
\noindent\textit{Proof.} 
Let $x,y\in G$ and, using Proposition \ref{brphaus}, take $N\in \N$ given by the BRP. Consider words $u=\overline{x^{-1}y}$ and $u'\in L_2$ such that $u'\pi_2=(x\varphi)^{-1}y\varphi$. We have that $$\Haus\left(\Img(\theta_u^x\varphi),\Img(\theta_{u'}^{x\varphi})\right)<N.$$ So, for every $k\in\{0,\ldots, |u'|\}$ we have that $$d_A(x\varphi u'^{[k]}\pi_2, \Img(\theta_u^x\varphi))\leq d_A(x\varphi u'^{[k]}\pi, G\varphi) <N.$$
\qed

We don't know if the converse of Proposition \ref{brphaus}  holds in general. We will now see that it does if $L_2$ is an f-biautomatic structure admitting a departure function.

\begin{Prop}
\label{brpstrong}
Let $G$ be a  biautomatic group with an automatic structure $L_1$ for $\pi_1:A^*\to G$ and an  f-biautomatic structure $L_2$ for $\pi_2:B^*\to G$ admitting some departure function $D$. Let $\varphi$ be an endomorphism such that there is some $N\in \N$ satisfying the following condition:
 \begin{itemize}[label=(C)]
 \item for every  $\alpha\in L_1$, $x\in G$, we have that $$ \Img(\theta_\alpha^x\varphi) \subseteq \mc V_{N}(\Img(\theta_{\beta}^{x\varphi}))<N,$$
for every $\beta$  such that $\beta\pi_2=\alpha\pi_1\varphi.$ 
\end{itemize}
\noindent Then the FT-BRP holds for $(\varphi,L_1,L_2)$.
\end{Prop}
\noindent\textit{Proof.}  Let $ B_\varphi=\max\{d_A(1,a\varphi)\mid a\in A\}$, $K$ given by Lemma \ref{kuniq} for $L_2$, $u_i \in L_1$ and $u'_i \in L_2 \cap u_i\pi_1\varphi\pi_2^{-1}$ for $i = 1,2,3$ and suppose that $u_1u_2,u_3$ are $p$-MFT in $\Gamma_A(G)$. In view of Proposition \ref{umetodos}, we may assume that $u_i'$ is a minimal length representative of  $u_i\pi_1\varphi$ in $L_2$. Then $d_A(u_1\pi_1,u_3^{[|u_1|]}\pi_1)<p$, and so, $d_A(u_1\pi_1\varphi,u_3^{[|u_1|]}\pi_1\varphi)<pB_\varphi$. By hypothesis, there is some $k\in \{0,\ldots, |u_3'|\}$ such that $d_A(u_3^{[|u_1|]}\pi_1\varphi,u_3'^{[k]}\pi_2)<N$ and so, 
$$d_B(u_1'\pi_2,u_3'^{[k]}\pi_2)\leq N_{A,B}d_A(u_1'\pi_2,u_3'^{[k]}\pi_2)=N_{A,B}d_A(u_1\pi_1\varphi,u_3'^{[k]}\pi_2)<N_{A,B}(N+pB_\varphi).$$ 
Consider the factorization of $u_3'$  given by 
$$\begin{tikzcd}
1 \ar[rr,"\gamma_1"] && u_3'^{[k]}\pi_2 \ar[rr,"\gamma_2"] && u_3'\pi_2.
\end{tikzcd}$$
Since $L_2$ is factorial, then $\gamma_1,\gamma_2\in L_2.$ So, $\gamma_1$ and $u_1'$ are words in $L_2$ with the same starting point that end at bounded distance. This means that there is some $p_1$ such that $\gamma_1$ and $u_1'$ $p_1$-fellow travel in $\Gamma_B(G)$. Similarly, $\gamma_2$ and $u_2'$ are two words in $L$ that start at bounded distance and have the same endpoint. Since $L_2$ is biautomatic, there is some $p_2$ such that $\gamma_2$ and $u_2'$ $p_2$-fellow travel.

By Lemma \ref{concatftp}, we have that $u_1'u_2'$ and $u_3'=\gamma_1\gamma_2$ fellow travel and the FT-BRP holds for $(\varphi,L_1,L_2).$ 
\qed

Observing that for a hyperbolic group we have that $Geo_A(G)$ is an f-biautomatic structure with a departure function, we have the following corollary that shows that the FT-BRP is equivalent to the formulations in Theorem \ref{equivsbrphiper}.
\begin{Cor}
Let $G$ be a hyperbolic group and $\varphi\in \End(G)$. Then the FT-BRP holds for $(\varphi,Geo_A(G),Geo_A(G))$ if and only if the hyperbolic BRP holds for $\varphi$. \qed
\end{Cor}

We now see that the definition of equivalent automatic structures can be relaxed, which will simplify some arguments.

\begin{Lem}
\label{asynchaus}
Two automatic structures $L_1$ and $L_2$ are equivalent if and only if $L_1\cup L_2$ has the Hausdorff closeness property.
\end{Lem}
\noindent\textit{Proof.} 
Let $L_1, L_2$ be automatic structures for $\pi_1:A^*\to G$ and $\pi_2:B^*\to G$. Clearly, if $L_1$ and $L_2$ are equivalent, then $L_1\cup L_2$ has the Hausdorff closeness property. Now, suppose that  $L_1\cup L_2$ has the Hausdorff closeness property. 
We want to prove that $L_1$ and $L_2$ are equivalent, i.e., $L=L_1\cup L_2\subseteq (A\cup B)^*$ is an asynchronous automatic structure for $\pi_3:(A\cup B)^*\to G$ defined naturally. Consider the boundedly asynchronous automatic structures $L_1'\subseteq L_1$ and $L_2'\subseteq L_2$ given by  \cite[Theorem 7.2.4]{[ECHLPT92]}. We have that $L_1'$ is equivalent to $L_1$ and $L_2'$ is equivalent to $L_2$, so we will prove that $L_1'$ is equivalent to $L_2'$ and that suffices by transitivity. To do so, we will prove that $L_1'\cup L_2'$ satisfies both conditions in Lemma \ref{async}.
Condition 2  is satisfied by hypothesis. 
We now prove the existence of a departure function for $L_1'\cup L_2'$. Let $p_1,p_2$ be the (synchronous) fellow traveller constants satisfied by $L_1$ and $L_2$ (and consequently by $L_1'$ and $L_2'$), respectively, and  $D_1,D_2$ be departure functions for the structures $(L_1',\pi_1)$, $(L_2',\pi_2)$, respectively. Take $D:\RR_0^+\to \RR_0^+$ defined by $$D(x)=\max\{D_1(N_{A,A\cup B}x),D_1(N_{B,A\cup B}x),D_2(N_{A,A\cup B}x),D_2(N_{B,A\cup B}x)\},$$ for every $x\in \RR_0^+$. Let $w\in L_1'\cup L_2'$, $r,s\geq 0$, $t\geq D(r)$ and $s+t\leq \lvert w\lvert $. Suppose w.l.o.g. that $w\in L_1'$. Then 
$$d_{A\cup B}(w^{[s]}\pi_3,w^{[s+t]}\pi_3)\geq \frac 1 {N_{A,A\cup B}}d_{A}(w^{[s]}\pi_1,w^{[s+t]}\pi_1)\geq r.$$ 
\qed

\section{Quasiconvex subgroups}
\label{secqci}
In this section, we will focus on endomorphisms whose image is $L$-quasiconvex in the sense of Gersten and Short (see \cite{[GS91]}). The main result of the section proves that, in some sense, the synchronous BRP holds for endomorphisms with finite kernel and $L$-quasiconvex image. 

As usual, for $n\in \N$, we will denote by $[n]$ the set $\{1,\ldots, n\}$.
Gersten and Short proved in \cite{[GS91]} that quasiconvex subgroups of automatic groups are automatic and gave a structure for the subgroup. We will replicate their definitions and follow the lines in their proof. Let $L\subseteq A^*$ be an automatic structure for $\pi_1:A^*\to G$ and $H$ be an $L$-quasiconvex subgroup with constant $k$. Take a  word $w=a_1\ldots a_n$ in $L'(H)=L\cap \pi_1^{-1}(H)$. Since $H$ is quasiconvex, and $w\pi_1\in H$, we have that for every $i\in [n]$ there is some word $g_i\in A^*$ such that $\lvert g_i\lvert \leq k$ and $w^{[i]}\pi_1 g_i\pi_1\in H$. So, $w\pi_1=\prod\limits_{i=1}^n (g_{i-1}\pi_1)^{-1}(a_i\pi_1) (g_i\pi_1)$, where $g_0=g_n=\varepsilon$. Hence, each element in $H$ can be written as a product of elements of norm at most $2k+1$. We take $B$ to be the set of those words together with their inverses and let $L''(H)$ be the set of words in $L'(H)$ rewritten as words in $B^*$. Notice that every $b\in B$ represents an element of $H$. When the quasiconvex subgroup $H$ is clearly set, we usually write $L'$ and $L''$ instead of $L'(H)$ and $L''(H)$, respectively. We remark that this notation will be adopted throughout the section, so, whenever we write $L'$ or $L''$, we will always be referring to this construction and will mostly be used when the subgroup is the image of an endomorphism.

We now present a technical lemma that will be very useful later on. 
\begin{Lem}
\label{distl''}
Let $G$ be an automatic group with automatic structure $L$ for $\pi:A^*\to G$ and $H$ be an $L$-quasiconvex subgroup with constant $k$ and automatic structure $L''$ for $\pi_2:B^*\to H$. There is $K\in \N$ such that, for all $x,y\in H$, $$\frac 1 {2k+1}d_A(x,y)\leq d_B(x,y)\leq Kd_A(x,y).$$
\end{Lem}
\noindent\textit{Proof.} Take $K\in \N$ given by Lemma  \ref{kuniq} applied to $L$ and let $x,y\in H$.
By construction of $L''$, elements of $B$ are words of length at most $2k+1$ in $A^*$. So, $$\frac 1 {2k+1}d_A(x,y)\leq d_B(x,y).$$ Take $w$ to be an arbitrary representative of $x^{-1}y$ of minimal length in $L$. By (\ref{minrep}), there is $K$ such that $\lvert w\lvert \leq Kd_A(x,y)$. By construction of the rewriting process, there is a word $w'\in L''\subseteq B^*$ of length equal to $\lvert w\lvert $ such that $w'\pi_2=x^{-1}y.$ Hence, $$d_B(x,y)\leq \lvert w'\lvert =\lvert w\lvert \leq Kd_A(x,y).$$
\qed

  \begin{Prop}
 \label{l'syncl''}
 Let $G$ be an automatic group with automatic structure $L$ for $\pi_1:A^*\to G$ and $H$ be an $L$-quasiconvex subgroup with constant $k$. Let $B$ be the canonical set of generators of $L''$ with $\pi_2:B^*\to H$ and $\pi_3:(A\cup B)^*\to G$, defined naturally. Then, there is some $N\in \N$ satisfying the following property:\\
 
 for words  $u,v\in L' \cup L''$ such that $d_{A}(u\pi_3,v\pi_3)\leq 1$, the inequality
 $$d_A(u^{[n]}\pi_3,v^{[n]}\pi_3)\leq N$$
 
holds for all $n\in \N$.
\end{Prop}
 
 \noindent\textit{Proof.} By \cite[Theorem 3.1]{[GS91]}, $L''$ is an automatic structure for $\pi_2:B^*\to H$ defined by $b\pi_2=b\pi_1$. Let $\tilde L=L'\cup L''$. Let $K$ be given by Lemma \ref{distl''}, $M$ and $M''$ be the constants given by the fellow traveller property of $L$ and $L''$, respectively, and take $u, v\in \tilde L$ such that $d_{A}(u\pi_3,v\pi_3)\leq1.$  If both $u$ and $v$ belong to $L'\subseteq L$ they $M$-fellow travel in $\Gamma_A(G)$ and if $u,v\in L''$, then $d_B(u\pi_2,v\pi_2)\leq K$ and so $u$ and $v$ $(KM'')$-fellow travel in $\Gamma_B(H)$. By Lemma \ref{distl''}, the result follows.
 
So suppose w.l.o.g. that $u\in L'$ and $v\in L''$. There is some word $w\in L''$ obtained by rewriting $u$ as a word in $B^*$. Also, by construction, we have that 
\begin{align}
\label{in1}
d_A(u^{[n]}\pi_1,w^{[n]}\pi_2)\leq k,
\end{align}
 for all $n\in \{0,\ldots \lvert u\lvert \}$, where $w$ is seen as a word in $B^*$. In this sense, $w^{[1]}$ is a word of length at most $2k+1$ in $A^*$. But then, $w$ is a word in $B^*$ such that $d_{A}(w\pi_3,v\pi_3)\leq 1$, thus $d_B(w\pi_2,v\pi_2)\leq K$. We know that paths starting in the same point ending at bounded distance fellow travel in $\Gamma_B(H)$. Hence, we have that
$$d_B(w^{[n]}\pi_2,v^{[n]}\pi_2)\leq KM''$$ and so  
\begin{align}
\label{in2}
d_{A}(w^{[n]}\pi_3,v^{[n]}\pi_3)\leq (2k+1)d_B(w^{[n]}\pi_2,v^{[n]}\pi_2)\leq (2k+1)KM''
\end{align}
for all $n\in \N$.
Combining (\ref{in1}) with (\ref{in2}), we get that 
$$d_{A}(u^{[n]}\pi_3,v^{[n]}\pi_3)=d_{A}(u^{[n]}\pi_1,v^{[n]}\pi_2)\leq (2k+1)KM''+k,$$
for all $n\in \N$. 
 \qed
 
 Let $G_1,G_2$ be automatic groups and $L$ be an automatic structure for $\pi:A^*\to G_1$. Consider $\varphi:G_1\to G_2$ to be a homomorphism. We say that  $L$ \emph{induces an automatic structure through $\varphi$} if $L$ is an automatic structure of $G\varphi$ for $\pi\varphi$. We call this the \emph{automatic structure induced by $L$ through $\varphi$} and denote it by $L^{(\varphi)}$.  In particular, the existence of an induced structure implies that $G\varphi$ is automatic. 
 
We recall that, given a group $G=\langle A\rangle$, the \emph{geodesic metric} $d_A$, defined by letting $d_A(x,y)$ to be the length of the shortest path in $\Gamma_A(G)$ connecting $x$ to $y$. In the proof of the following Lemma we will always specify how letters in $A$ are read to avoid confusion, since we are dealing with the same language read through different surjective homomorphisms.

\begin{Lem}
\label{quotient}
Let $G$ be an automatic group, $L$ be an automatic structure for $\pi_1:A^*\to G$ and $\varphi\in \End(G)$ be a virtually injective endomorphism of $G$. Then $\faktor G {\Ker(\varphi)}$ is automatic and $L$ is an automatic structure for $\pi_1\pi$, where $\pi:G\to \faktor G {\Ker(\varphi)}$ denotes the projection onto the quotient.
\end{Lem}
\noindent\textit{Proof.} Let $K$ denote $\Ker(\varphi)$. We have that $G=\langle A\pi_1 \rangle$ and $\faktor G K=\langle A\pi_1\pi \rangle$. Take $M=\max\{d_{A\pi_1}(1,x)\mid x\in K\}$. Since $\pi_1$ and $\pi$ are sujective, then so is $\pi_1\pi$. We will now see that the fellow traveller property holds. So, take two words  $u,v\in L$ such that $d_{A\pi_1\pi}(u\pi_1\pi,v\pi_1\pi)\leq 1$. Then, there is some $a\in A\cup\{1\}$ such that $u\pi_1\pi=(v\pi_1\pi)(a\pi_1\pi)$, which means, by definition of $\pi$, that $u\pi_1\varphi=(va)\pi_1\varphi$ and so $(u\pi_1)^{-1}(va)\pi_1\in K$. By definition of $M$, $d_{A\pi_1}(u\pi_1,(va)\pi_1)=d_{A\pi_1}(1,(u\pi_1)^{-1}(va)\pi_1)\leq M$. So, there is some $N$ given by the fellow traveller property for $\pi_1$ for $L$-geodesics ending at distance at most $M$ such that 
$$d_{A\pi_1}(u^{[n]}\pi_1, v^{[n]}\pi_1) \leq N$$
holds for all $n\in \N$
and so 
$$d_{A\pi_1\pi}(u^{[n]}\pi_1\pi, v^{[n]}\pi_1\pi) \leq N$$
holds for all $n\in \N$.
\qed

\begin{Cor}
Let $G$ be an automatic group and $\varphi\in \End(G)$ be virtually injective. Then $L$ induces an automatic structure through $\varphi$. 
\end{Cor}
\noindent\textit{Proof.} Let $\varphi':G\to G\varphi$ be the homomorphism obtained by restricting the codomain of $\varphi$ to the image. Then $\varphi'=\pi\bar\varphi$, where $\bar\varphi:\faktor{G}{\Ker(\varphi)}\to G\varphi$ is an isomorphism and so taking an automatic structure $L$ for $\pi_1:A^*\to G$, by Lemma \ref{quotient}, we have that $L$ is an automatic structure for $\pi_1\pi\bar\varphi$.
\qed

We remark that the converse does not hold. To see that, it suffices to consider an endomorphism with finite image.

\begin{Thm}
\label{brp}
Let $G_1$ and $G_2$ be automatic groups with automatic structures $L_1$ and $L_2$ for $\pi_1:A^*\to G_1$ and $\pi_2:B^*\to G_2$, respectively. Let $\varphi:G_1\to G_2$ be a homomorphism such that $G_1\varphi$ is $L_2$-quasiconvex and suppose that $L_1$ induces an automatic structure through $\varphi.$ Then $L_1^{(\varphi)}$ and $L_2''(G_1\varphi)$ are (synchronous) equivalent if and only if the (synchronous) BRP holds for $(\varphi, L_1,L_2).$
\end{Thm}

\noindent\textit{Proof.}  Let $C$ be the standard  set of generators of $L_2''$.

\textit{Asynchronous case.} Suppose that $L_1^{(\varphi)}$ and $L_2''$ are equivalent and put $L=L_1^{(\varphi)}\cup L_2''$. Then $L$ is an asynchronous automatic structure for $\pi_3:(A\cup C)^*\to G_1\varphi$. 
Take 
$p$ to be the asynchronous fellow travel constant satisfied by $L$ and $q$ to be the quasiconvexity constant. 

Let $\alpha\in L_1,$ $x\in G_1$, and $\beta\in L_2$  be such that $\beta\pi_2=\alpha\pi_1\varphi$, $k\in \{0,\ldots,\lvert \alpha\lvert \}$ and  $y=(\theta_\alpha^x(k))\varphi$.  Take $N\in \N$ given by Proposition \ref{l'syncl''} and $\beta'\in L_2''$ such that $\beta'\pi_3=\alpha\pi_1\varphi$. 
By Proposition \ref{l'syncl''}, we have that  $$d_B(\beta^{[n]}\pi_2,\beta'^{[n]}\pi_3)\leq N$$
holds for all $n\in \N$.
By hypothesis, we have that $\alpha$, when read in $L_1^{(\varphi)}$, asynchronously $p$-fellow travels $\beta'$ in $\Gamma_{A\cup C}(G_1\varphi)$. 
Let $\phi,\psi$ be the reparametrisation functions for $\alpha$ and $\beta'$ and take $t\in \N$ such that $\phi(t)=k$. This way, we have that $y=(x\alpha^{[\phi(t)]\pi_1})\varphi$. By construction, elements of $C$ are words of length at most $2q+1$ in $B^*$. Hence, putting $\lambda= \max(\{2q+1\}\cup \{d_B(1,a\pi_1\varphi) \mid a\in A\})$, we have that 
\begin{align*}
d_B\left(y,\Img(\theta_{\beta}^{x\varphi})\right) 
\leq \;  &d_B\left(y,x\varphi\left(\beta^{[\psi(t)]}\pi_2\right)\right)\\
 \leq \; & d_B\left(x\varphi\left(\alpha^{[\phi(t)]}\pi_1\varphi\right),x\varphi\left(\beta'^{[\psi(t)]}\pi_2\right)\right) \\
&+   d_B\left(x\varphi\left(\beta'^{[\psi(t)]}\pi_2\right),x\varphi\left(\beta^{[\psi(t)]}\pi_2\right)\right)\\
\leq 
\;&\lambda d_{A\cup C}\left(\alpha^{[\phi(t)]}\pi_1\varphi,\beta'^{[\psi(t)]}\pi_2\right) + N\\
\leq \;&\lambda p + N
\end{align*}
Now, let $k'\in \{0,\ldots \lvert \beta\lvert \}$, $z=\theta_{\beta}^{x\varphi}(k')$ and $k''$ such that $\psi(k'')=k'$. This way, $z=x\varphi(\beta^{[\psi(k'')]}\pi_2)$. We have that 

\begin{align*}
d_B\left(z,\Img(\theta_\alpha^x\varphi)\right)\leq \; & d_B\left(z,x\varphi\left(\alpha^{[\phi(k'')]}\pi_1\varphi\right)\right)\\
 \leq \; & d_B\left(z,x\varphi\left(\beta'^{[\psi(k'')]}\pi_2\right)\right) +   d_B\left(x\varphi\left(\beta'^{[\psi(k'')]}\pi_2\right),x\varphi\left(\alpha^{[\phi(k'')]}\pi_1\varphi\right)\right)\\
\leq  \; & d_B\left(x\varphi\left(\beta^{[\psi(k'')]}\pi_2\right),x\varphi\left(\beta'^{[\psi(k'')]}\pi_2\right)\right)  \\
&+\lambda  d_{A\cup C}\left(x\varphi\left(\beta'^{[\psi(k'')]}\pi_2\right),x\varphi\left(\alpha^{[\phi(k'')]}\pi_1\varphi\right)\right) \\
\leq \; & N +\lambda p
\end{align*}

It follows that $$\Haus\left(\Img(\theta_\alpha^x\varphi),\Img(\theta_{\beta}^{x\varphi})\right)\leq N + \lambda p+1.$$
%%%%%%%%%%%
%%CONVERSE%%%
%%%%%%%%%%%
Now, suppose that $L_1\subseteq A^*$ and $L_2\subseteq B^*$ are such that there is some $K\in \N$ satisfying the following condition: for all $\alpha\in L_1$, $x\in G$, we have that $$\Haus\left(\Img(\theta_\alpha^x\varphi),\Img(\theta_{\beta}^{x\varphi})\right)\leq K,$$
for every $\beta\in L_2$  such that $\beta\pi_2=\alpha\pi_1\varphi.$
We want to check that $L_1^{(\varphi)}$ and $L_2''$ are equivalent. In view of Lemma \ref{asynchaus}, it suffices to see that $L_1^{(\varphi)}\cup L_2''$ has the Hausdorff closeness property. 

Consider two words $w_1,w_2\in L_1^{(\varphi)}\cup L_2'',$ whose images under the map $\pi_3$ are at distance at most one apart in $\Gamma_{A\cup C}(G_1\varphi)$.

If $w_1,w_2\in L_1^{(\varphi)}$ or $w_1,w_2\in L_2''$, then their Hausdorff distance is bounded by the synchronous fellow traveller constant of $L_1^{(\varphi)}$ or $L_2''$, respectively. So, suppose w.l.o.g. that $w_1\in L_1^{(\varphi)}$ and $w_2\in L_2''.$ 
Take $w_2''\in L_2''$ such that $w_2''\pi_2=w_1\pi_1\varphi$ and $w_2'\in L_2'\subseteq L_2$ such that $w_2''$ is obtained by rewriting $w_2'$ in $C^*$. By hypothesis, we have that $$\Haus\left(\Img(\theta_{w_1}^1\varphi),\Img(\theta_{w_2'}^{1})\right)\leq K.$$  
Also, by construction, $d_B(w_2''^{[n]}\pi_2,w_2'^{[n]}\pi_3)\leq q$, for all $n\in \N$, thus $w_1$ and $w_2''$ are $(K+q)$-Hausdorff close in $\Gamma_B(G_2)$. By Lemma \ref{distl''}, it follow that they are Hausdorff close in $\Gamma_{C}(G_1\varphi)$. Since $w_2''$ and $w_2$ are two words in $L_2''$ ending at distance at most one apart, they fellow travel in $\Gamma_C(G_1\varphi)$. So, $w_1$ and $w_2$ are Hausdorff close in $\Gamma_{C}(G_1\varphi)$ and thus in $\Gamma_{A\cup C}(G_1\varphi).$  

Hence, $L_1^{(\varphi)}$ and $L_2''$ are equivalent.\\

 \textit{Synchronous case.} Suppose that $L_1^{(\varphi)}$ and $L_2''$ are synchronous equivalent and put $L=L_1^{(\varphi)}\cup L_2''$. Then $L$ is an automatic structure for $\pi_3:(A\cup C)^*\to G_1\varphi$. 
Take 
$p$ to be the synchronous fellow travel constant satisfied by $L$ and $q$ to be the quasiconvexity constant. 

Let $\alpha\in L_1$ and $\beta\in L_2$  be such that $\beta\pi_2=\alpha\pi_1\varphi$ and $n\in\N$.  We want to prove that $d_B(\alpha^{[n]}\pi_1\varphi,\beta^{[n]}\pi_2)$ is bounded.

Take $N\in \N$ given by Proposition \ref{l'syncl''} and $\beta'\in L_2''$ obtained by rewriting $\beta$ in $C^*$. Then, by construction, we have that
$$d_B(\beta'^{[n]}\pi_3,\beta^{[n]}\pi_2)\leq q.$$

By hypothesis, we have that $\alpha$, when read in $L_1^{(\varphi)}$, synchronously $p$-fellow travels $\beta'$ in $\Gamma_{A\cup C}(G_1\varphi)$, so $$d_{A\cup C}(\alpha^{[n]}\pi_1\varphi, \beta'^{[n]}\pi_3)\leq p.$$
Hence, putting $\lambda= \max(\{2q+1\}\cup \{d_B(1,a\pi_1\varphi) \mid a\in A\})$, we have that 
$$d_{B}(\alpha^{[n]}\pi_1\varphi, \beta'^{[n]}\pi_3)\leq \lambda p.$$
Thus, $$d_B(\alpha^{[n]}\pi_1\varphi,\beta^{[n]}\pi_2)\leq q+\lambda p.$$

%%%%%%%%%%%
%%CONVERSE%%%
%%%%%%%%%%%
To prove the converse, suppose that the synchronous BRP holds for $(\varphi,L_1,L_2)$ with constant $K$. We want to prove that $L=L_1^{(\varphi)}\cup L_2''$ is an automatic structure for $G_1\varphi$ with the homomorphism $\pi_3:(A\cup C)^*\to G_1\varphi$ defined in the natural way. Take words $\alpha, \beta\in L$ such that $d_{A\cup C}(\alpha\pi_3,\beta\pi_3)\leq 1$. If both $\alpha$ and $\beta$ belong to either $L_1^{(\varphi)}$ or $L_2''$, then they fellow travel in $\Gamma_A(G_1\varphi)$ or $\Gamma_C(G_1\varphi)$, respectively and so they do in $\Gamma_{A\cup C}(G_1\varphi)$. So suppose w.l.o.g. that $\alpha\in L_1^{(\varphi)}$ and $\beta\in L_2''$.
Take $\gamma''\in L_2''$ such that $\gamma''\pi_3=\alpha\pi_3$ and $\gamma'\in L_2'\subseteq L_2$ such that $\gamma''$ is obtained by rewriting $\gamma'$ in $C^*$.
 For every $n\in \N$, we have that $$d_B(\alpha^{[n]}\pi_3,\gamma'^{[n]}\pi_2)=d_B(\alpha^{[n]}\pi_1\varphi,\gamma'^{[n]}\pi_2)\leq K,$$ by hypothesis. Also, by construction, we have that 
 $$d_B(\gamma'^{[n]}\pi_2, \gamma''^{[n]}\pi_3)\leq  q,$$ for all $n\in \N$. Applying \ref{difgen}, we have that $$d_{A\cup C}(\alpha^{[n]}\pi_3, \gamma''^{[n]}\pi_3)\leq d_A(\alpha^{[n]}\pi_3, \gamma''^{[n]}\pi_3)\leq  N_{A,B}(q+K)$$ and since $\gamma''$ and $\beta$ are words in $L_2''$ ending at bounded distance, they synchronously fellow travel in $\Gamma_C(G_1\varphi)$ and so, they do in $\Gamma_{A\cup C}(G_1\varphi)$ and the result follows. 
\qed

\begin{Cor}
\label{idbrp}
Let $G$ be an automatic group with automatic structures $L_1$ and $L_2$ for $\pi_1:A^*\to G$ and $\pi_2:B^*\to G$ respectively, and consider the identity mapping $\iota:(G,L_1)\to (G,L_2)$. Then  $L_1$ and $L_2$ are (synchronous) equivalent if and only if the (synchronous) BRP holds for $(\iota, L_1,L_2).$ In particular, if $G$ is hyperbolic, then the BRP holds for $(\iota,L_1,L_2)$.
\end{Cor}

\begin{Cor}
Let $G$ be an automatic group with automatic structure $L$ for $\pi_1:A^*\to G$ and $\varphi\in \End(G)$ be an endomorphism inducing an automatic structure $L^{(\varphi)}$ on $G\varphi$. Then the (synchronous) BRP holds for $(\varphi',L,L^{(\varphi)})$, where $\varphi':G\to G\varphi$ is the homomorphism obtained by restricting the codomain of $\varphi$ to the image.
\end{Cor}

Combining Corollary \ref{idbrp} with Lemma \ref{compos}, we get that, if $G$ is hyperbolic, then the BRP is independent of the structures considered, since we can compose with the identity mapping on the left (resp. right) to change the structure considered in the domain (resp. codomain) and the BRP will be preserved. Also, if $G$ is an automatic group with automatic structure $L$ for $\pi_1:A^*\to G$ and $\varphi\in \End(G)$ is an endomorphism with $L$-quasiconvex image such that $L$ induces an automatic structure through $\varphi$,  then the BRP holds for $(\varphi,L,L)$ if and only if  $L^{(\varphi)}$ equivalent to $L''$. 
In particular, if the image is hyperbolic then the BRP must hold for  $(\varphi,L,L)$. This yields the already known result that if $G$ is a hyperbolic group and $\varphi\in \End(G)$ is virtually injective and has quasiconvex image, i.e., a quasi-isometric embedding (see \cite[Theorem 4.3]{[AS16]}) then the BRP holds for $\varphi$.  We will prove something more general later in this section.

We now prove a weaker result for the FT-BRP.

\begin{Prop}
\label{strongbrp}
Let $G$ be an automatic group with automatic structures $L_1$ and $L_2$ for $\pi_1:A^*\to G$ and $\pi_2:B^*\to G$, respectively. Let $\varphi\in \End(G)$ be such that $G\varphi$ is $L_2$-quasiconvex with constant $q$ and $L_1$ induces an automatic structure through $\varphi.$ If $L_1^{(\varphi)}$ and $L_2''$ are synchronous equivalent and $L_1^{(\varphi)}\cup L_2''$ admits a departure function, then the FT-BRP holds for $(\varphi, L_1,L_2).$
\end{Prop}
\noindent\textit{Proof.}  Let $C$ be the standard  set of generators of $L_2''$ and $\pi''$ be the surjective homomorphism $\pi'':C^*\to G\varphi$. Put $L=L_1^{(\varphi)}\cup L_2''$. We have that $L$ is a synchronous automatic structure for $\pi_3:(A\cup C)^*\to G\varphi$. Take 
$s$ to be the fellow travel constant satisfied by $L$. 
Also, let $N$ be the constant given by Proposition \ref{l'syncl''}, M be given by Lemma \ref{concatftp} with $p_1=p_2=s$, and $p\geq 0$.

Let $u_i \in L_1$ and $u'_i \in L_2 \cap u_i\pi_1\varphi\pi_2^{-1}$ for $i = 1,2,3$ such that  $u_1u_2,u_3$ are $p$-MFT in $\Gamma_A(G)$. Clearly, $\{u_1',u_2',u_3'\}\subseteq L_2'.$  For $i=1,2,3$, let $u_i''$ be a word in  $L_2''$ obtained by rewriting $u_i'$. Then, by Proposition \ref{l'syncl''}, we have that, for all $i=1,2,3,$ $n\in \N$, 
$$d_A(u_i'^{[n]}\pi_2,u_i''^{[n]}\pi'')<N.$$

 By hypothesis, $u_i''$ and $u_i$ $s$-fellow travel in $\Gamma_{A\cup C}(G\varphi)$ when $u_i$ is read by $\pi_1\varphi.$ Also, $u_i'\in B^*$ and $u_i''\in C^*$ have the same length by construction.

Hence, we have that 
\begin{align}
\label{ii1}
d_B((u_1'u_2')^{[n]}\pi_2,(u_1''u_2'')^{[n]}\pi'' )<N.
\end{align}
Since $u_i''$ and $u_i$ $s$-fellow travel in $\Gamma_{A\cup C}(G\varphi)$ and $L_1^{(\varphi)}\cup L_2''$ admits a departure function, then by Lemma \ref{concatftp}, $u_1''u_2''$ and $u_1u_2$ $M$-fellow travel in $\Gamma_{A\cup C}(G\varphi)$.  
Since elements in $C$ are words  in $B^*$ of size at most $2q+1$ and elements in $A$ when read through $\pi_1\varphi$ can be represented by words in $B^*$ of size at most $\max\{d_B(1,a\varphi)\mid a\in A\}$, putting $$\lambda=\max\{2q+1,\max\{d_B(1,a\varphi)\mid a\in A\}\},$$ we have that 
\begin{align}
\label{ii2}
d_B((u_1''u_2'')^{[n]}\pi'',(u_1u_2)^{[n]}\pi_1\varphi)<\lambda M.
\end{align}
 So, combining (\ref{ii1}) and (\ref{ii2}), we get that
\begin{align}
\label{ii3}
d_B((u_1'u_2')^{[n]}\pi_2,(u_1u_2)^{[n]}\pi_1\varphi)<N+\lambda M.
\end{align}
 Now, let $B_\varphi=\max\{d_A(1,a\varphi)\mid a\in A\}$. We have that $u_1u_2$ and $u_3$ $p$-fellow travel in $\Gamma_A(G)$ and so
$$d_A((u_1u_2)^{[n]}\pi_1\varphi,u_3^{[n]}\pi_1\varphi)<pB_\varphi$$
for all $n\in \N$, thus
\begin{align}
\label{ii4}
d_B((u_1u_2)^{[n]}\pi_1\varphi,u_3^{[n]}\pi_1\varphi)<N_{A,B}pB_\varphi.
\end{align}
 But $u_3$ and $u_3''$ $s$-fellow travel in $\Gamma_{A\cup C}(G\varphi)$ and so 
\begin{align}
\label{ii5}
d_B(u_3''^{[n]}\pi'',u_3^{[n]}\pi_1\varphi)<s\lambda
\end{align}
and by Proposition  \ref{l'syncl''} 
\begin{align}
\label{ii6}
d_B(u_3''^{[n]}\pi'',u_3'^{[n]}\pi_2)<N.
\end{align}

Combining (\ref{ii3}), (\ref{ii4}), (\ref{ii5}) and (\ref{ii6}), we get that  $u'_1u'_2,u'_3$ are MFT in $\Gamma_B(G)$.
\qed

In particular, if two structures $L_1$ and $L_2$ are synchronous equivalent and admit departure functions, then the FT-BRP holds for the identity mapping $\iota:(G,L_1)\to (G,L_2)$. 

\begin{Cor}
Let $G$ be an automatic group with automatic structure $L$ for $\pi_1:A^*\to G$ and $\varphi\in \End(G)$ be an endomorphism inducing an automatic structure $L^{(\varphi)}$ on $G\varphi$ admitting a departure function. Then the FT-BRP holds for $(\varphi',L,L^{(\varphi)})$, where $\varphi':G\to G\varphi$ is the endomorphism obtained by restricting the codomain of $\varphi$ to the image.
\end{Cor}
\noindent\textit{Proof.} Let $p>0$, $B_\varphi= \max\{d_A(1,a\varphi)\mid a\in A\}$ and take $u_i \in L_1$  such that $u_1u_2,u_3$ are $p$-MFT in $\Gamma_A(G)$. Then  $u_i \in L_1^{(\varphi)}\cap u_i\pi\varphi(\pi\varphi)^{-1}$ for $i = 1,2,3$ and obviously  $u_1u_2,u_3$ are $pB_\varphi$-MFT in $\Gamma_A(G\varphi)$. By Proposition \ref{umetodos}, the result follows.
\qed

Using Theorem \ref{brp} and Proposition \ref{strongbrp}, we can find a connection between the BRP and the FT-BRP for a larger class of structures, generalizing Propositions \ref{brphaus} and \ref{brpstrong}.

\begin{Prop}
Let $G$ be an automatic group,  $L_1$ be an automatic structure for $\pi_1:A^*\to G$ admitting a departure function, $L_2$ be an automatic structure for $\pi_2:B^*\to G$ and $\varphi\in \End(G)$. Then, if the FT-BRP holds for $(\varphi,L_1,L_2)$, so does the BRP.
\end{Prop}
\noindent\textit{Proof.} Consider the identity mapping $\iota$. By Proposition \ref{strongbrp}, the FT-BRP holds for $(\iota, \text{Fact($L_1$}), L_1)$. Since it holds for $(\varphi,L_1,L_2)$, then it holds for $(\varphi,\text{Fact($L_1$}), L_2)$  by Lemma \ref{compos}. By Proposition \ref{brphaus}, the BRP holds for $(\varphi,\text{Fact($L_1$}), L_2)$. Since $L_1$ and Fact($L_1$) are equivalent, then the BRP holds for $(\iota,L_1,\text{Fact($L_1$)})$ by Theorem \ref{brp}. Thus, it holds for $(\varphi,L_1,L_2)$ by Lemma \ref{compos}.
\qed

\begin{Prop}
Let $G$ be a biautomatic group,  $L_1$ be an automatic structure for $\pi_1:A^*\to G$, $L_2$ be a biautomatic structure for $\pi_2:B^*\to G$ %synchronous equivalent to $L_1$, 
admitting a departure function,
and $\varphi\in \End(G)$. Then, if the BRP holds for $(\varphi,L_1,L_2)$, so does the FT-BRP.
\end{Prop}
\noindent\textit{Proof.} Consider the identity mapping $\iota$. Since Fact($L_2$) is (synchronous) equivalent to $L_2$ and it also admits a departure function, then the BRP holds for $(\iota,L_2,\text{Fact($L_2$)})$ and the FT-BRP holds for $(\iota,\text{Fact($L_2$)},L_2)$.
 Using Lemma \ref{compos}, we get that the BRP holds for $(\varphi, L_1, \text{Fact($L_2$)})$. By Proposition \ref{brpstrong}, the FT-BRP holds for $(\varphi, L_1, \text{Fact($L_2$)})$. By Lemma \ref{compos}, it holds for $(\varphi, L_1, L_2).$
\qed\\

Finally, we can prove that in some sense, the synchronous BRP always holds for virtually injective endomorphisms with quasiconvex image.

\begin{Thm}
\label{estsync}
Let $G$ be an automatic group with an automatic stucture $L$ for $\pi_1:A^*\to G$ and $\varphi$ be a virtually injective endomorphism with $L$-quasiconvex image. Then there is some automatic structure $\tilde L$ such that the synchronous BRP holds for $(\varphi, \tilde L, L)$.
\end{Thm}

\noindent\textit{Proof.}
Let $q$ be the quasiconvexity constant, $K$ be given by Lemma \ref{distl''} and put $V_\varphi=\max\{d_A(1,a\varphi)\mid a\in A\}$ and $M=\max\{d_A(1,g)\mid g\in \Ker(\varphi)\}$. Consider $L''=L''(G\varphi)$, which is an automatic structure for $\pi_2:B^*\to G\varphi$ \cite[Theorem 3.1]{[GS91]}. Put $$S=B\pi_2=\{b\pi_2\mid b\in B\}.$$ We fix a total ordering of $A$. For every $s\in S$, $s\varphi^{-1}$ is finite and we denote the shortlex minimal word in $A^*$ that represents an element in $s\varphi^{-1}$ by $\tilde s_{\varphi^{-1}}$. Let $$C=\{\tilde s_{\varphi^{-1}}\mid s\in S\},$$
which is finite since $S$ is finite. Notice that it might be the case where $\varepsilon\in S$. Also, let $N_C=\max\{d_A(1,c\pi_1)\mid c\in C\}$ and $$L_\Ker=\{\overline g \in A^*\mid g\in \Ker(\varphi)\},$$ which is obviously finite.  Now, we will define a language $\tilde L$ on $(A\cup C)^*$ and prove that it defines an automatic structure for $\pi_3:(A\cup C)^*\to G$ defined by $w\pi_3=w\pi_1$, where letters in $C$ are viewed as words in $A^*$, such that $\tilde L^{(\varphi)}$ is synchronous equivalent to $L''$. Take a word $w=w_1w_2\cdots w_n\in L''$ and put $x_i=w_i\pi_2$, for $i=1,\ldots, n.$ Consider all words of the form $\widetilde {x_1}_{\varphi^{-1}}\cdots \widetilde {x_n}_{\varphi^{-1}}w\in (A\cup C)^*$, where $w\in L_\Ker.$ So, each word is rewritten in $\lvert \Ker(\varphi)\lvert $ different ways and let $\tilde L$ be the language of the words obtained by rewriting all words in $L''$. We will  prove that: 
\begin{enumerate}
\item $\tilde L$ is rational;
\item $\pi_3\lvert _{\tilde L}$ is a surjective homomorphism;
\item  $\tilde L$ satisfies the fellow traveler property;
\item $\tilde L^{(\varphi)}$ and $L''$ are synchronous equivalent.
\end{enumerate}

To prove $1$, let $\mc A''=(Q, q_0, T, \delta)$ be a finite state automaton recognizing $L''$. Then we replace every transition labelled by $b\in B$ by a transition labelled $\widetilde{(b\pi_2)}_{\varphi^{-1}}$. Then we add a new terminal state $q_T$ and for each terminal state $q\in T$, add paths from $q$ to $q_T$ labelled by all words in $L_\Ker$. The language accepted by the new automaton is precisely $\tilde L$.

To prove $2$, let $g\in G$ and take a word $w=w_1\cdots w_n\in L''$ representing $g\varphi$. Put $x_i=w_i\pi_2$, for $i=1,\ldots, n$ and consider  the word  $w'=\widetilde {x_1}_{\varphi^{-1}}\cdots \widetilde {x_n}_{\varphi^{-1}}\in C^*$. Then, by definition of $\widetilde {x_i}_{\varphi^{-1}}$ and since $\varphi$ is a homomorphism, we have that $w'$ represents (via $\pi_3$) an element $g'\in G$ such that $g'\varphi=g\varphi$, i.e., there is some $k\in \Ker(\varphi)$ such that $g=g'k$. By construction of $\tilde L$, $w'\bar k\in \tilde L$ and it represents $g$ when read through $\pi_3$. Since $g$ is arbitrary, we have that $\pi_3\lvert _{\tilde L}$ is a surjective homomorphism.

To prove $3$, take two words $u,v\in \tilde L$ such that $d_{A\cup C}(u\pi_3,v\pi_3)\leq 1$. Then, by construction $u=u_Cu_A$ and $v=v_Cv_A$ where $u_C,v_C\in C^*$ and $u_A,v_A\in L_\Ker$.

Let $u'',v''\in L''$ be words from which $u,v$ are obtained through rewriting. 
\begin{align*}
d_B(u''\pi_2,v''\pi_2)&=d_B(u\pi_3\varphi,v\pi_3\varphi)\leq Kd_A(u\pi_3\varphi,v\pi_3\varphi)\\
&\leq KV_\varphi d_A(u\pi_3,v\pi_3) \\
&\leq KV_\varphi N_C d_{A\cup C}(u\pi_3,v\pi_3) \\
&\leq KV_\varphi N_C.
\end{align*}

Since $L''$ is an automatic structure, then there is some $N$ depending only on $K,V_\varphi, N_C$ such that $u''=u_1\cdots u_r$ and $v''=v_1\cdots v_s$ $N$-fellow travel in $\Gamma_B(G\varphi)$. Take $$P=\max\{d_A(1,g)\mid g\in G, d_A(1,g\varphi)\leq (2q+1)N\},$$ which is well defined since the kernel is finite. We now claim that $u_C$ and $v_C$ $(P+M)$-fellow travel in $\Gamma_{A\cup C}(G)$. Indeed, let $n\in \N$. For $i>r$, put $u_i=\varepsilon$ and for $i>s$, put $v_i=\varepsilon$. Then,
by Lemma \ref{distl''}, we have that 
\begin{align*}
&d_{A}\left(\left(\widetilde{(u_1\pi_2)}_{\varphi^{-1}}\cdots \widetilde{(u_n\pi_2)}_{\varphi^{-1}}\right)\pi_3\varphi,\left(\widetilde{(v_1\pi_2)}_{\varphi^{-1}}\cdots \widetilde{(v_n\pi_2)}_{\varphi^{-1}}\right)\pi_3\varphi\right)\\
\leq &(2q+1)d_{B}\left(\left(\widetilde{(u_1\pi_2)}_{\varphi^{-1}}\cdots \widetilde{(u_n\pi_2)}_{\varphi^{-1}}\right)\pi_3\varphi,\left(\widetilde{(v_1\pi_2)}_{\varphi^{-1}}\cdots \widetilde{(v_n\pi_2)}_{\varphi^{-1}}\right)\pi_3\varphi\right)\\
= &(2q+1)d_{B}(u''^{[n]}\pi_2,v''^{[n]}\pi_2)\\
\leq & (2q+1)N.
\end{align*}
Thus, there are letters $a_1,\ldots, a_k\in A$ for some $k\leq (2q+1)N$ such that 
$$\left(\widetilde{(u_1\pi_2)}_{\varphi^{-1}}\cdots \widetilde{(u_n\pi_2)}_{\varphi^{-1}}\right)\pi_3\varphi (a_1\cdots a_k)\pi_3= \left(\widetilde{(v_1\pi_2)}_{\varphi^{-1}}\cdots \widetilde{(v_n\pi_2)}_{\varphi^{-1}}\right)\pi_3\varphi$$
and $(a_1\cdots a_k)\pi_3\in G\varphi$. Since $k\leq (2q+1)N$, then  $(a_1\cdots a_k)\pi_3=\alpha\pi_3\varphi$, for some $\alpha\in A^*$ such that 
\begin{align}
\label{lenalpha}
\lvert \alpha\lvert \leq P.
\end{align} 

Hence, 
$$
\left(\left(\widetilde{(u_1\pi_2)}_{\varphi^{-1}}\cdots \widetilde{(u_n\pi_2)}_{\varphi^{-1}}\alpha\right)\pi_3\right)^{-1}\left(\widetilde{(v_1\pi_2)}_{\varphi^{-1}}\cdots \widetilde{(v_n\pi_2)}_{\varphi^{-1}}\right)\pi_3\in \Ker(\varphi),$$ so 
\begin{align}
\label{lenker}
d_A\left(\left(\widetilde{(u_1\pi_2)}_{\varphi^{-1}}\cdots \widetilde{(u_n\pi_2)}_{\varphi^{-1}}\alpha\right)\pi_3,\left(\widetilde{(v_1\pi_2)}_{\varphi^{-1}}\cdots \widetilde{(v_n\pi_2)}_{\varphi^{-1}}\right)\pi_3\right)\leq M 
\end{align}

Using (\ref{lenalpha}) and (\ref{lenker}), we get that
\begin{align*}
&d_A\left(\left(\widetilde{(u_1\pi_2)}_{\varphi^{-1}}\cdots \widetilde{(u_n\pi_2)}_{\varphi^{-1}}\right)\pi_3,\left(\widetilde{(v_1\pi_2)}_{\varphi^{-1}}\cdots \widetilde{(v_n\pi_2)}_{\varphi^{-1}}\right)\pi_3\right)\\
\leq \;&  d_A\left(\left(\widetilde{(u_1\pi_2)}_{\varphi^{-1}}\cdots \widetilde{(u_n\pi_2)}_{\varphi^{-1}}\right)\pi_3,\left(\widetilde{(u_1\pi_2)}_{\varphi^{-1}}\cdots \widetilde{(u_n\pi_2)}_{\varphi^{-1}}\alpha\right)\pi_3\right)\\
&+ d_A\left(\left(\widetilde{(u_1\pi_2)}_{\varphi^{-1}}\cdots \widetilde{(u_n\pi_2)}_{\varphi^{-1}}\alpha\right)\pi_3,\left(\widetilde{(v_1\pi_2)}_{\varphi^{-1}}\cdots \widetilde{(v_n\pi_2)}_{\varphi^{-1}}\right)\pi_3\right)\\
\leq \;  & d_A(1,\alpha\pi_3)+M\\
\leq\; &\lvert \alpha\lvert +M\\
\leq\; & P+M.
\end{align*}
Therefore, 
\begin{align*}
&d_{A\cup C}(u_C^{[n]}, v_C^{[n]})\\
=\;&d_{A\cup C}\left(\left(\widetilde{(u_1\pi_2)}_{\varphi^{-1}}\cdots \widetilde{(u_n\pi_2)}_{\varphi^{-1}}\right)\pi_3,\left(\widetilde{(v_1\pi_2)}_{\varphi^{-1}}\cdots \widetilde{(v_n\pi_2)}_{\varphi^{-1}}\right)\pi_3\right)\\
\leq \; & d_{A}\left(\left(\widetilde{(u_1\pi_2)}_{\varphi^{-1}}\cdots \widetilde{(u_n\pi_2)}_{\varphi^{-1}}\right)\pi_3,\left(\widetilde{(v_1\pi_2)}_{\varphi^{-1}}\cdots \widetilde{(v_n\pi_2)}_{\varphi^{-1}}\right)\pi_3\right)\\
\leq  \; &P+M.
\end{align*}

It follows that $u$ and $v$ $(P+3M)$-fellow travel because $\lvert u_A\lvert ,\lvert v_A\lvert \leq M$.

Now, $4$ is essentially obvious by construction. Consider $\tilde L^{(\varphi)}\cup L''$ and define $\pi_4:(A\cup C\cup B)^*\to G\varphi$ naturally. Take words $u,v\in \tilde L^{(\varphi)}\cup L''$ such that $d_{A\cup C\cup B}(u\pi_4,v\pi_4)\leq 1$. If they both belong to $\tilde L$ or $L''$, then they fellow travel in $\Gamma_{A\cup C}(G\varphi)$ or $\Gamma_B(G\varphi),$ respectively and so, they do in  $\Gamma_{A\cup C\cup B}(G\varphi)$. So, suppose w.l.o.g that $u\in \tilde L$ and $v\in L''$ and consider the factorization $u=u_Cu_A$ as done above. Then $u\pi_4=u\pi_3\varphi=u_C\pi_4$ since $u_A\pi_4\in \Ker(\varphi)$. Consider the word $u''\in L''$ from which $u$ is obtained through rewriting. Since $u''$ and $v$ fellow travel in $\Gamma_B(G\varphi)$, so do $u_C$ (when read through $\pi_3\varphi$) and $v$. Since $\lvert u_A\lvert \leq M$, then $u$ and $v$ fellow travel in  $\Gamma_B(G\varphi)$, and so they do in  $\Gamma_{A\cup C\cup B}(G\varphi)$.

These four points, combined with Theorem \ref{brp}, yield the desired result.
\qed

\section{Some applications}
The goal of this section is to apply the techniques developed in the previous sections in order to prove quasiconvexity of interesting subgroups defined by endomorphisms.

\label{secapps}
Let $A,B$ be finite alphabets. Take the alphabet $$C=(A\times B)\cup (A\times\{\$\})\cup (\{\$\}\cup B)$$ and define the convolution of two words $u\in A^*, v\in B^*$ as the only word in $C^*$ whose projection on the first (resp. second) component belongs to $u\$^*$ (resp. $v\$^*$).  The convolution can be defined naturally for languages $K\subseteq A^*$ and $L\subseteq B^*$ by taking $$K\diamond L=\{u\diamond v\mid u\in K, \,v\in L\}.$$

We start by defining a natural structure on the direct product of two automatic groups.

\begin{Prop}
\label{convprod}
Let $G_1$ and $G_2$ be automatic groups and take automatic structures $L_1$ for $\pi_1:A^*\to G_1$ and $L_2$ for $\pi_2:B^*\to G_2$. 
Then $L_1\diamond L_2$ is an automatic structure for $G_1\times G_2$.
\end{Prop}
\noindent\textit{Proof.} We start by proving that $L_1\diamond L_2$ is a rational language of $C^*$. Define the obvious projection-like homomorphisms $\rho_1:C^*\to A^*$ and $\rho_2:C^*\to B^*$. We have that $$A^*\diamond B^*=(A\times B)^*((A\times \{\$\})^*\cup (\{\$\}\times B)^*)$$ and so it is rational. Since $$K\diamond L=K\rho_1^{-1}\cap L\rho_2^{-1}\cap (A^*\diamond B^*),$$ it is also rational.
Now, define $\pi_3:C^*\to G_1\times G_2$ by $(x,y)\mapsto (x\pi_1,y\pi_2)$ if $(x,y)\in A\times B$, $(x,\$)\mapsto (x\pi_1,1)$,and $(\$,y)\mapsto (1,y\pi_2)$. It is clearly surjective and $\pi_3\lvert _{G_1\diamond G_2}$ is still surjective. We only have to check that the fellow traveler property holds. So, take words $x=u_1\diamond v_1$ and $y=u_2\diamond v_2$ in $L_1\diamond L_2$ such that $d_C(x\pi_3,y\pi_3)\leq 1.$
Notice that $d_C(x\pi_3,y\pi_3)=\max\{d_A(u_1\pi_1,u_2\pi_1),d_B(v_1\pi_2,v_2\pi_2)\}$ and so there is some $r>0$  such that $u_1,u_2$ $r$-fellow travel in $\Gamma_A(G_1)$ and $v_1,v_2$ $r$-fellow travel in $\Gamma_B(G_2)$. Let $n\in \N$. We have that $$d_C(x^{[n]}\pi_3,y^{[n]}\pi_3)=\max\{d_A(u_1^{[n]}\pi_1,u_2^{[n]}\pi_1),d_B(v_1^{[n]}\pi_2,v_2^{[n]}\pi_2)\}\leq r.$$
\qed

\begin{Prop}
\label{kerfg}
Let $G$ be an automatic group with automatic structures $L_1$ and $L_2$ for $\pi_1:A^*\to G$ and $\pi_2:B^*\to G$, respectively.
Consider an endomorphism  $\varphi:G\to G$ such that $L_1$ induces an automatic structure on $G\varphi$ and the BRP holds for $(\varphi, L_1,L_2)$. Then, $\Ker(\varphi)$ is isomorphic to a $(L_1\diamond L_1^{(\varphi)})$-quasiconvex  subgroup of $G\times G\varphi$ . In particular, $\Ker(\varphi)$ is automatic.
\end{Prop}
\noindent\textit{Proof.} We start by showing that we can assume that $L_2\cap1\pi_2^{-1}=\{\varepsilon\}$. If the BRP holds for $(\varphi, L_1,L_2)$, then it holds for $(\varphi, L_1,L_2')$, where $L_2'\subseteq L_2$ is an automatic structure with uniqueness. If we replace the unique representative of $1$ in $L_2'$ by $\varepsilon$, we obtain a new automatic structure $L_3$ equivalent to $L_2$, and so, by Corollary \ref{idbrp}, the BRP holds for $(\varphi, L_1,L_3)$.

Let $N\in \N$ be the constant given by the BRP. Since $L_1$ induces an automatic structure on $G\varphi$, then by Proposition \ref{convprod}, we have that $G\times G\varphi$ is an automatic group and $L_1\diamond L_1^{(\varphi)}$ is an automatic structure of $G\times G\varphi$. Now, put $H=\{(x,x\varphi)\mid x\in G\}\leq G\times G\varphi$. We have that $$\Ker(\varphi)\simeq \{(x,x\varphi)\mid x\in \Ker(\varphi)\}=(G\times \{1\}) \cap H.$$
Now we will prove that both $H$ and $G\times\{1\}$ are $(L_1\diamond L_1^{(\varphi)})$-quasiconvex subgroups of $G\times G\varphi$. 
It is obvious by construction that $H$ is $(L_1\diamond L_1^{(\varphi)})$-quasiconvex. Indeed, take a word $u\diamond v\in L_1\diamond L_1^{(\varphi)}$ representing $(x,x\varphi)$ for some $x\in G$. We have that $u\pi_1\varphi=x\varphi=v\pi_1\varphi$ and so $u$ and $v$ fellow travel in $\Gamma_A(G\varphi)$, reading through $\pi_1\varphi$. Thus, $$d\left((u^{[n]}\pi_1,v^{[n]}\pi_1\varphi), H\right)\leq d\left((u^{[n]}\pi_1,v^{[n]}\pi_1\varphi),(u^{[n]}\pi_1,u^{[n]}\pi_1\varphi)\right),$$ which is bounded by the fellow traveler constant satisfied by $L_1^{(\varphi)}.$

Also, to prove that $G\times \{1\}$ is 
 $(L_1\diamond L_1^{(\varphi)})$-quasiconvex with constant $N$, observe that, taking a word $u\diamond v\in L_1\diamond L_1^{(\varphi)}$ representing an element in $G\times\{1\}$, we have that $v$ must be a word in $L_1$ such that $v\pi_1\in \Ker(\varphi).$ The path defined by $v$ in $\Gamma_A(G\varphi)$ where letters are read through $\pi_1\varphi$ is at a Hausdorff distance smaller than $N$ from $\{1\}$, i.e., $d_B(1,v^{[n]}\pi_1\varphi,1)\leq N$, because $\varepsilon$ is the only representative of $1$ in $L_2$ and the BRP holds for $(\varphi,L_1,L_2)$.
 Thus, every point of the path defined by $u\diamond v \in L_1\diamond L_1^{(\varphi)}$ is of the form $(g,h)$, where $g\in G$ and $d_B(1,h)\leq N$. % since the path defined by $v$ is at Hausdorff distance at most $N$ from $\{1\}$ since the BRP holds for $(\varphi, L_1,L_2)$ with constant $N$. 
 Letting $M=\max \{d_A(1,x)\mid x\in G\varphi ,\, d_B(1,x) \leq N\}$, where letters of $A$ are read through $\pi_1\varphi$, we have that $G\times\{1\}$ is $(L_1\diamond L_1^{(\varphi)})$-quasiconvex with quasiconvexity $M$. 
 
 Since both subgroups are $(L_1\diamond L_1^{(\varphi)})$-quasiconvex, they are $(L_1\diamond L_1^{(\varphi)})$-rational and so their intersection is also $(L_1\diamond L_1^{(\varphi)})$-rational, thus $(L_1\diamond L_1^{(\varphi)})$-quasiconvex.
\qed

\begin{Prop}
Let $G$ be an automatic group with automatic structures $L_1$ and $L_2$ for $\pi_1:A^*\to G$ and $\pi_2:B^*\to G$, respectively  and consider an endomorphism $\varphi\in \End(G)$ such that the synchronous BRP holds for $(\varphi, L_1,L_2)$. Then $L_1$ induces an automatic structure through $\varphi$.
\end{Prop}
\noindent\textit{Proof.} Let $N$ be the constant given by the synchronous BRP for $ (\varphi, L_1,L_2)$, $V_\varphi=\max\{d_B(1,a\varphi)\mid a\in A\}$ and $K$ be the fellow traveller constant satisfied by words in $L_2$ ending at distance at most $V_\varphi$. The only thing we need to check is that $L_1^{(\varphi)}$ satisfies the fellow traveller property. So let $u,v \in L_1$ be such that $d_A(u\pi_1\varphi, v\pi_1\varphi)\leq 1$. Let $u',v'\in L_2$ be such that $u'\pi_2=u\pi_1\varphi$ and $v'\pi_2=v\pi_1\varphi$. 
Since the synchronous BRP holds for $(\varphi, L_1,L_2)$, we have that, for all $n\in \N$,  $$d_B(u^{[n]}\pi_1\varphi,u'^{[n]}\pi_2) \leq N \quad \text{ and } d_B(v^{[n]}\pi_1\varphi,v'^{[n]}\pi_2)\leq N.$$
Also, we have $d_B(u'\pi_2,v'\pi_2)\leq V_\varphi$, thus $u'$ and $v'$ $K$-fellow travel in $\Gamma_B(G)$. So, for every $n\in \N$, we have that 
\begin{align*}
d_B(u^{[n]}\pi_1\varphi,v^{[n]}\pi_1\varphi)\leq \; &d_B(u^{[n]}\pi_1\varphi,u'^{[n]}\pi_2)+  d_B(u'^{[n]}\pi_2,v'^{[n]}\pi_2)\\
&+ d_B(v'^{[n]}\pi_2,v^{[n]}\pi_1\varphi)\\
 \leq \; & 2N+K
\end{align*}
and so, letting $M=\max \{d_A(1,x)\mid x\in G\varphi ,\, d_B(1,x)\leq 2N+K\}$, where letters in $A$ are read through $\pi_1\varphi$, we have that $L_1^{(\varphi)}$ satisfies the $M$-fellow traveler property.
\qed

\begin{Rmk} 
\label{syncfgker}
The last two propositions combined show that having the synchronous BRP for any pair of languages $L_1$ and $L_2$ 
is enough to have an automatic, thus finitely presented, kernel. Indeed, if the synchronous BRP holds for $(\varphi,L_1,L_2)$, then it holds for $(\varphi,L_1,L_2')$ where $L_2'$ is an automatic structure with uniqueness by Corollary \ref{idbrp}. If we replace the unique representative of $1$ in $L_2'$ by $\varepsilon$, we obtain a new automatic structure $L_3$ such that the synchronous BRP holds for $(\varphi, L_1,L_3)$. Now, $\varphi$, $L_1$ and $L_3$ satisfy the hypothesis of Proposition \ref{kerfg}. This immediately shows that this is a very strong condition to impose on a general endomorphism. For example, in the case of virtually free groups, having finitely generated kernel implies that either the kernel or the image is finite. 
\end{Rmk}
We will prove something similar for fixed points with the additional hypothesis that $L_1$ and $L_2$ belong to the same synchronous equivalence class, which seems a strong condition to impose.

\begin{Prop}
\label{theta}
Let $G_1$ and $G_2$ be automatic groups with automatic structures $L_1$ and $L_2$ for $\pi_1:A^*\to G_1$ and $\pi_2:B^*\to G_2$, respectively. Let $\varphi,\psi:G_1\to G_2$ be  homomorphisms such that the synchronous BRP holds for $(\varphi, L_1,L_2)$ and $(\psi, L_1,L_2)$. Let   $\theta:G_1\to G_2\times G_2$ be the homomorphism defined by $x\mapsto (x\varphi,x\psi).$ Then the synchronous BRP holds for $(\theta,L_1,L_2\diamond L_2)$.
\end{Prop}
\noindent\textit{Proof.} Let $N_1$ and $N_2$ be the constants given by the  synchronous BRP holding for $(\varphi, L_1,L_2)$ and $(\psi, L_1,L_2)$, respectively, take $N=\max\{N_1,N_2\}$, and put, as usual, $C=(A\times B)\cup (A\times\{\$\})\cup (\{\$\}\cup B)$ and  $\pi_3:C^*\to G\times G$. Let $x\in G_1$,  $u\in L_1$ representing $x$ and $v\diamond w\in L_2\diamond L_2$ representing $(x\varphi,x\psi)$. Since the synchronous BRP holds for $(\varphi, L_1,L_2)$ and $(\psi, L_1,L_2)$, we have that, for all $n\in \N$,  $$d_B(u^{[n]}\pi_1\varphi, v^{[n]}\pi_2)\leq N$$ and  $$d_B(u^{[n]}\pi_1\psi, w^{[n]}\pi_2)\leq N.$$
Thus, $$d(u^{[n]}\pi_1\theta, (v\diamond w)^{[n]}\pi_3)=\max\{d_B(u^{[n]}\pi_1\varphi, v^{[n]}\pi_2),d_B(u^{[n]}\pi_1\psi, w^{[n]}\pi_2)\}\leq N$$
and the synchronous BRP holds for $(\theta,L_1,L_2\diamond L_2)$.
\qed

\begin{Ex} The fact that the BRP is synchronous for both homomorphisms is crucial to the proof above. Indeed, consider $\ZZ$ with the structure $L$ given by the geodesics and let $\theta:\ZZ\to \ZZ\times \ZZ$ defined by $a\mapsto (a,a^2)$. Then, the image is not $(L\diamond L)$-quasiconvex (and so the BRP does not hold for $(\varphi, L,L\diamond L)$) despite being of the form $x\mapsto (x\varphi,x\psi)$, where the synchronous BRP holds for $(\varphi,L,L)$ and the BRP holds for $(\psi,L,L)$. 
\end{Ex}

\begin{Cor}
\label{eq}
Let $G_1$ and $G_2$ be automatic groups with automatic structures $L_1$ and $L_2$ for $\pi_1:A^*\to G_1$ and $\pi_2:B^*\to G_2$, respectively. Let $\varphi,\psi:G_1\to G_2$ be  homomorphisms such that the synchronous BRP holds for $(\varphi, L_1,L_2)$ and $(\psi, L_1,L_2)$. 
Then $\text{Eq}(\varphi,\psi)=\{x\in G_1\mid x\varphi=x\psi\}$ is isomorphic to a $(L_2\diamond L_2)$-quasiconvex subgroup of $G_2\times G_2$. In particular, $\text{Eq}(\varphi,\psi)$ is automatic.
\end{Cor}
\noindent\textit{Proof.} 
It is obvious that the diagonal subgroup $\Delta=\{(x,x)\in G_2\times G_2: x\in G_2\}$ is an $(L_2\diamond L_2)$-quasiconvex subgroup of $G_2\times G_2$. The subgroup $H=\{(x\varphi,x\psi)\mid x\in G_1\}$ is also  $(L_2\diamond L_2)$-quasiconvex since it is the image of the homomorphism in Proposition \ref{theta}, for which the (synchronous) BRP holds. 
Since $\text{Eq}(\varphi)\simeq \Delta\cap H$, the result follows.
\qed

\begin{Cor}
\label{fix}
Let $G$ be an automatic group with synchronous equivalent automatic structures $L_1$ and $L_2$ for $\pi_1:A^*\to G$ and $\pi_2:B^*\to G$, respectively. Let $\varphi$ be an endomorphism such that the synchronous BRP holds for $(\varphi, L_1,L_2)$. 
Then $\Fix(\varphi)$  is isomorphic to a $(L_2\diamond L_2)$-quasiconvex subgroup of $G\times G$. In particular, $\text{Fix}(\varphi)$ is automatic.
\qed
\end{Cor}

In view of Corollary \ref{idbrp}, the hypothesis needed to apply Corollary \ref{fix}  is equivalent to the existence of some automatic structure $L$ such that the synchronous BRP holds for $(\varphi,L,L)$. This seems to be quite strong since it requires the endomorphism to distort points only a bounded amount, in some sense. We will now prove that for biautomatic groups, inner automorphisms satisfy this property for any biautomatic structure and see what happens in case the group is free and $L$ is the structure given by the geodesics.

\begin{Prop}
\label{innersynch}
Let $G$ be a biautomatic group, $L$ be a biautomatic structure for $\pi:A^*\to G$ and $\varphi\in \Inn(G)$ be an inner automorphism of $G$. Then, the synchronous BRP holds for $(\varphi, L,L)$.
\end{Prop}

\noindent\textit{Proof.}
For $a\in A$, let $\lambda_a$ be the inner automorphism given by $x\mapsto axa^{-1}$. We will see that, for every $a\in A$, the synchronous BRP holds for $(\lambda_a,L,L)$ and the result will follow by Lemma \ref{compos}.  Let $a\in A$, $w\in L$, $N$ be the fellow traveler property for paths in $\Gamma_A(G)$ labelled by words in $L$ starting and ending at distance at most one and take $w'\in L$ such that $w'\pi=a(w\pi)a^{-1}$. Consider the paths $\alpha$ starting in $1$ labelled by $w$ and $\beta$ starting in $a^{-1}$ labelled by $w'$. They start and end at distance at most $1$ (it is exactly one if $a\pi\neq 1 $), thus they $N$-fellow travel in $\Gamma_A(G).$ This means that for all $n\in \N$, $$d_A\left(a^{-1}(w'^{[n]}\pi), w^{[n]}\pi\right)\leq N,$$ 
and so 
$$d_A\left(w'^{[n]}\pi, a(w^{[n]}\pi)\right) \leq N,$$
then having
$$d_A\left(w'^{[n]}\pi, w^{[n]}\pi\varphi\right)=d_A\left(w'^{[n]}\pi, a(w^{[n]}\pi)a^{-1}\right)\leq N+1.$$
\qed

Even though the hypothesis of Corollary \ref{fix} are strong, it yields an alternative proof to Proposition 4.3 in \cite{[GS91]}.

\begin{Cor}
\label{centralizer}
 The centralizer of a finite subset of a biautomatic group is biautomatic.
\end{Cor}
\noindent\textit{Proof.} The centralizer of an element is the fixed subgroup of the inner automorphism defined by that element. Since, by Proposition \ref{innersynch}, inner automorphisms of biautomatic groups satisfy the synchronous BRP for $(\varphi,L,L)$, where $L$ is a biautomatic structure, then we can apply Corollary \ref{fix} to  get that the centralizer is quasiconvex (and thus, biautomatic). Now, the centralizer of a finite subset is a finite intersection of quasiconvex subgroups, and so quasiconvex.
\qed

We now prove that in the case of free groups, not many endomorphisms besides the inner automorphisms satisfy the hypothesis of Corollary \ref{fix}.

Given a free group, we call \emph{letter permutation automorphism} to an automorphism that maps generators into generators bijectively.

\begin{Prop}
 For a free group $F=F_A$, the endomorphisms $\varphi$ for which the synchronous BRP holds for $(\varphi, Geo_A(F), Geo_A(F))$ are precisely the automorphisms in the subgroup generated by the inner and the letter permutation automorphisms.
\end{Prop}
\noindent\textit{Proof.} Let $A$ be a finite alphabet, $F=F_A$ be the free group over $A$ and let $H$ be the subgroup of $\Aut(F)$   generated by the inner and the letter permutation automorphisms. For $u\in F_n$, we will denote by $\lambda_u$  the inner automorphism defined by $x\mapsto uxu^{-1}$.

It is obvious that the synchronous BRP holds for automorphisms consisting of permutations of letters with constant $0$ and Proposition \ref{innersynch} shows that it also holds for inner automorphisms, since $Geo_A(F)$ is a biautomatic structure of $F$.

So, by Lemma \ref{compos}, if $\varphi$ belongs to $H$, then the synchronous BRP holds for $(\varphi,Geo_A(F),Geo_A(F))$.

Now, we will prove the converse. Let $\varphi\in \End(F)$ be an endomorphism and let $a\in A$. Suppose that the cyclically reduced core of $a\varphi$ has length greater than $1$.  We have that $\lvert a^{n}\varphi\lvert -n>n$, and so, letting $w$ be a reduced word representing $a^n\varphi$, it follows that $d_A(w^{[n]},a^n\varphi)$ is unbounded, thus the synchronous BRP cannot hold for $(\varphi, Geo_A(F), Geo_A(F))$. So, if the synchronous BRP holds for $(\varphi, Geo_A(F), Geo_A(F))$, then the image of a letter $x$ is of the form $w_xy_xw_x^{-1}$ for some $w_x\in F$, $y_x\in A$. 

If, for every $x\in A$, $w_x$ is trivial, then the endomorphism is induced by a permutation of the generators. Notice that, in this case, we must have that $y_x\neq y_z$ for all $x,z\in A$ with $x\neq z$, because otherwise we would have that $(x^nz^{-n})\varphi=1$ and $d((x^nz^{-n})^{[n]}\varphi,1)=d(x^n\varphi,1)=d(y_x^n,1)=n$, which is unbounded and so the synchronous BRP cannot hold for $(\varphi, Geo_A(F), Geo_A(F))$. 

So suppose that there is some $x\in A$ such that $w_x\neq \varepsilon$.
Let $z\in A$ be such that $\lvert w_z\lvert=\max\{\lvert w_a\lvert\,\big\lvert\, a\in A\}$.  If there is some $x\in A$ such that $w_x\neq \varepsilon$ and $w_x$ is not a prefix of  $w_z$, then the cyclically reduced core of $(zx)\varphi=w_zy_zw_z^{-1}w_xy_xw_x^{-1}$ has length greater than $3$. Thus,
\begin{align*}
\lvert (zx)^n\varphi\lvert =\lvert (w_zy_zw_z^{-1}w_xy_xw_x^{-1})^n\lvert > 3n
\end{align*}
and so, letting $w$ be a geodesic representing $(zx)^n\varphi$, we have that $d_A((zx)^n\varphi,w^{[2n]})$ is unbounded. So, if the synchronous BRP holds for $(\varphi, Geo_A(F), Geo_A(F))$, then for every $x\in A$, $w_x$ is a prefix of $w_z$ (it might be the case where $w_x=\varepsilon$). 

We now proceed by induction on $m=\max\{\lvert w_a\lvert\,\big\lvert\, a\in A\}$. If $m=0$, then we are done.
Suppose now that for an endomorphism $\varphi$ such that  $m \leq n$, if the synchronous BRP  holds for $(\varphi, Geo_A(F), Geo_A(F))$, then $\varphi$ belongs to $H$. Take $\varphi$ such that  $\max\{\lvert w_a\lvert\,\big\lvert\, a\in A\}=n+1$, fix $z\in A$ such that $\lvert w_z \lvert=n+1$ and let $b$ be the first letter in $w_z$.  

Suppose that there is some $x\in A$ such that $w_x=\varepsilon$. In this case, if $y_x\neq b^{-1},b$, then $(zx)\varphi=w_zy_zw_z^{-1}y_x$ is cyclically irreducible and has length greater than $3$, thus the same argument as above shows that the synchronous BRP cannot hold for $\varphi$. This means that every $x\in A$ having $w_x=\varepsilon$ must be such that $y_x=b^{-1}$ or $y_x=b.$ Obviously, we have that $\varphi=(\varphi\lambda_{b}^{-1})\lambda_b$. The image of a letter $x$ through $\varphi\lambda_{b}^{-1}$ is equal to $v_xy_xv_x^{-1}$, where $v_x$ is just $w_x$ without the first letter (might be empty), if $\lvert w_x\lvert\neq \varepsilon$, and it is equal to $b$ or $b^{-1}$ if $x\varphi=b$ or $x\varphi=b^{-1}$, respectively. In view of Lemma \ref{compos} and Proposition \ref{innersynch}, the synchronous BRP holds for $\varphi$ if and only if it holds for $\varphi\lambda_{b}^{-1}$. Applying the induction hypothesis to $\varphi\lambda_{b}^{-1}$, we get that if the synchronous BRP holds for $\varphi\lambda_{b}^{-1}$, then $\varphi\lambda_{b}^{-1}\in H$ and so does $\varphi$.
\qed

 We don't know if the same also holds for structures not synchronous equivalent to $Geo_A(F)$. Also, answering the same question for larger classes of groups is yet to be done.
\begin{Question}
Given an automatic structure $L$ on a free group $F$, is it true that the synchronous BRP holds for $(\varphi,L,L)$ if and only if $\varphi$ is in the subgroup generated by the inner automorphisms and the ones induced by permutations of the generators? 
\end{Question}

\begin{Question}
Given an automatic (biautomatic, hyperbolic, virtually free) group, which endomorphisms $\varphi$ are such that there exists an automatic structure $L$ such that the synchronous BRP holds for $(\varphi,L,L)$?
\end{Question}

In \cite{[GS91]}, Gersten and Short give an example of an automorphism of a biautomatic group whose fixed subgroup is not finitely generated: letting $G = F_2\times \ZZ$, the automorphism $\varphi$ of $G$ given by $(a,0)\mapsto (a,0)$, $b\mapsto (b,1)$ and $(0,1)\mapsto (0,1)$ has fixed subgroup $N\times \ZZ$, where $N$ is the normal closure of $a$ in $F_2$ . Thus the fixed subgroup of $\varphi$ is not finitely generated. In Theorem \ref{estsync}, it is proved that, given an automorphism $\varphi$, there are structures $L_1$ and $L_2$ such that the synchronous BRP holds for $(\varphi, L_1, L_2)$. However, $L_1$ and $L_2$ might not be synchronous equivalent, as this example shows.

\section{Classification of endomorphisms}
\label{secclass}
In \cite{[Car21]}, we considered the properties of having quasiconvex image, satisfying the BRP, being uniformly continuous for a visual metric and being virtually injective and obtained implications between them and gave counterexamples for the ones that do not hold for nontrivial endomorphisms of hyperbolic groups. As seen above, in hyperbolic groups, quasiconvexity and the BRP are independent of the automatic structures. When the synchronous BRP is considered or the class of groups is expanded to the whole class of automatic groups, this is not the case. So, from now on, when we say that the (synchronous) BRP holds for an endomorphism, we mean that it holds for some pair of structures $L_1$ and $L_2$ and having quasiconvex image will mean that there is some automatic structure $L$ such that the image of the endomorphism is $L$-quasiconvex.. We will consider in detail the cases of free, virtually free, hyperbolic and automatic groups.

\subsection{Free groups}

In free groups, every endomorphism has quasiconvex image, and being injective is equivalent to the BRP. From \ref{estsync}, it follows that it also coincides with the synchronous BRP. 

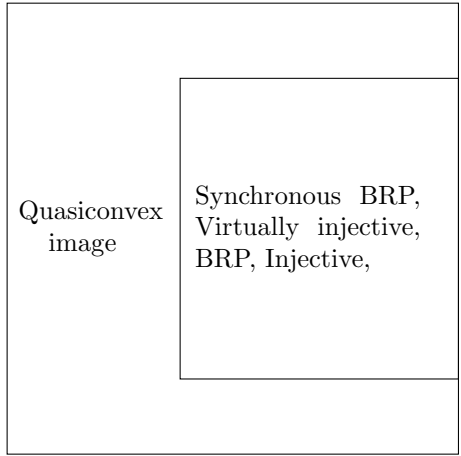
\begin{figure}[H]
\begin{center}
\begin{tikzpicture}
	\begin{scope} [fill opacity = 0.4]

    \draw[ draw = black] (-4.5,-2) rectangle (1.5,4);
        \draw[ draw = black] (-2.2,-1) rectangle (1.5,3);

   \node[opacity=1] at (-3.5,1) {\small\text{\parbox{1.7cm} {\centering Quasiconvex image}}};
    \node[opacity=1] at (-0.5,1) {\small\text{\parbox{3cm} {Synchronous BRP, Virtually injective, BRP,  Injective, }}};

    \end{scope}

\end{tikzpicture}
\caption{Nontrivial endomorphisms of free groups}
\end{center}
\end{figure}

\subsection{Virtually free groups}

We know that the BRP holds for every virtually injective endomorphism of a virtually free group (and so the image is quasiconvex). From \ref{estsync}, it follows that the synchronous BRP also holds. The synchronous BRP holds for endomorphisms with finite image. We will now prove that is the only possible case of endomorphisms with infinite kernel for which the  
BRP holds.

Let  $G$ be a virtually free group and consider a decomposition as a disjoint union
$$ G= Fb_0\cup Fb_1\cup \cdots \cup Fb_m,$$
where $F_A\trianglelefteq G$ is a free group of finite rank and $1=b_0,b_1,\ldots, b_m\in G.$ Let $\beta_0,\ldots, \beta_m$ be letters outside of $A$ and put $B=\{\beta_0,\ldots,\beta_m\}$ and $C=A\cup A^{-1}\cup B$. We start by proving that $L=Geo_A(F)B\subseteq C^*$ is an automatic structure for $\pi:C^*\to G$ defined by $x\mapsto x$ if $x\in A\cup A^{-1}$ and $\beta_i\mapsto b_i$, for $\in\{0,\ldots, m\}$. 

For $i\in \{0,\ldots, m\}$ and $g\in G$, let $\sigma(i,g)$ be the element of $F$ such that $b_ig=\sigma(i,g)b_k$ for some $k\in[m]$  and put $$M=\max\{|\sigma(i,w)|\,\lvert\, i\in \{0,\ldots,m\},\, w\in A\cup A^{-1}\cup\{b_0,\ldots, b_m\} \}.$$ Since both  $Geo_A(F)$ and $B$ are regular sublanguages of $C^*$, then so is their concatenation. Also, $\pi$ is a surjective homomorphism, so we only have to verify that the fellow traveler property holds for $L$. Consider two words $u\beta_i, v\beta_j\in L$ such that $d_C((u\beta_i)\pi,(v\beta_j)\pi)\leq 1$. Then, there is some $x\in C$ such that $ub_i=(u\beta_i)\pi=(v\beta_jx)\pi$. This means that $\sigma(j,x)=v^{-1}u$, and so $$d_A(u\pi,v\pi)\leq |\sigma(j,x)|\leq M,$$
and so $u$ and $v$ fellow travel in $\Gamma_A(F)$. Thus, we have that $d_C((u\beta_i)^{[n]}\pi,(v\beta_j)^{[n]}\pi)\leq d_A(u^{[n]}\pi,v^{[n]}\pi)$ if $n\leq \min\{|u|,|v|\}$ and, if not, by the triangular inequality, it follows that  $d_C((u\beta_i)^{[n]}\pi,(v\beta_j)^{[n]}\pi)\leq d_A(u^{[n]}\pi,v^{[n]}\pi)+2$ and so $u\beta_i$ and $v\beta_j$ fellow travel in $\Gamma_C(G)$.

Now, let $\varphi\in \End(G)$ be an endomorphism which has infinite kernel and infinite image. Then, $\varphi|_F$ is neither trivial nor injective.  Then $F\varphi$ must contain an element of infinite order. Indeed, since $F\varphi$ is infinite, there must be some nontrivial $w\in F\cap F\varphi$  since there must be some $i\in [m]$ and $g,h\in F$ such that $g\neq h$ and $gb_i,hb_i\in F\varphi$ and so $gh^{-1}\in F\varphi\setminus\{1\}$. Thus, there are nontrivial $u\in F$ such that $u\varphi=1$ and $v\in F$ such that $v\varphi$ has infinite order. We may assume that $u$ is cyclically reduced (if it is not we replace $u$ by its cyclically reduced core).
%, that $u$ is not a power of $v$, and that $u$ is neither a prefix nor a suffix of $v$. Since $u$ is cyclically reduced,
Therefore, for every $n\in \N$ there cannot be cancellation boyh in $v^{n}u$ and in $v^{n}u^{-1}$. 

 Suppose that the BRP holds with constant $N$. Take $n\in \N$ such that $d_C(v^n\varphi, 1)>N$ and $d_C(v^{-n}\varphi, 1)>N$. Suppose w.l.o.g. that cancellation occurs in  $v^{n}u$ and take $k\in \N$ such that $|u^k|>|v^n|$. This way, $v^n$ is a prefix of $v^nu^{-k}v^{-n}$.
Hence, we have that $(v^nu^{-k}v^{-n})\varphi=1$, but we have that $d_C(1,v^n\varphi)>N$, which contradicts the BRP. So, if the BRP holds for $\varphi$ (since the group is hyperbolic the BRP is not dependent on the structures we take), then either the kernel or the image must be finite.

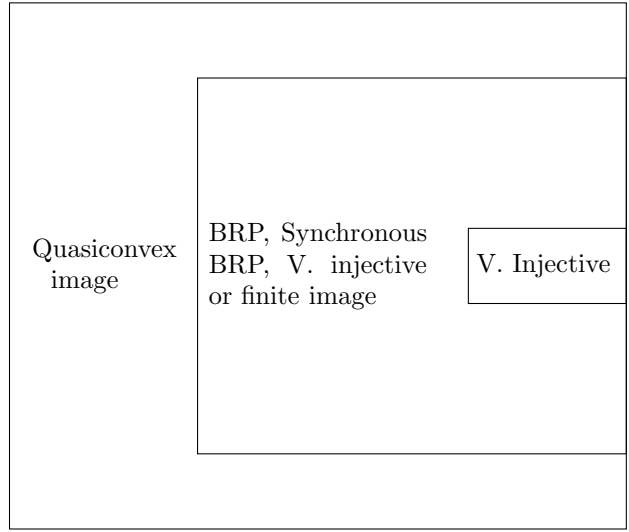
\begin{figure}[H]
\begin{center}
\begin{tikzpicture}
	\begin{scope} [fill opacity = 0.4]

    \draw[ draw = black] (-6.5,-2.5) rectangle (1.7,4.5);
        \draw[ draw = black] (-4,-1.5) rectangle (1.7,3.5);
\draw[ draw = black] (-0.4,0.5) rectangle (1.7,1.5);

   \node[opacity=1] at (-5.5,1) {\small\text{\parbox{1.4cm} {\centering Quasiconvex image}}};

    \node[opacity=1] at (-2.4,1) {\small\text{\parbox{2.9cm} {BRP, Synchronous BRP, V. injective or finite image}}};

        \node[opacity=1] at (0.6,1) {\small\text{V. Injective}};
    \end{scope}

\end{tikzpicture}
\caption{Nontrivial endomorphisms of virtually free groups}
\end{center}
\end{figure}

\subsection{Hyperbolic groups}

In the case of hyperbolic groups, again, the synchronous BRP must hold for virtually injective endomorphism with quasiconvex image 
by Theorem \ref{estsync}. However, unlike the virtually free groups case, Example 5.4 in \cite{[AS16]} shows that there is a virtually injective endomorphism of a torsion-free hyperbolic group with a non quasiconvex image.
 The same kind of questions as the ones in the virtually free groups arise.

\begin{Question}
\label{q1}
Is there an endomorphism of a hyperbolic group for which the BRP holds and the synchronous BRP does not? 
\end{Question} 
\begin{Question}
Is there an endomorphism of a hyperbolic group for which the BRP holds but the kernel is not finitely generated? By Remark \ref{syncfgker}, an affirmative answer to this question would also yield an affirmative answer to Question \ref{q1}.
\end{Question} 

In virtually free groups, we know that finitely generated normal subgroups are finite or have finite index, so having a finitely generated kernel is the same as having a finite kernel or a finite image. We wonder if something similar to the virtually free groups case
might hold in the case of hyperbolic groups.

\begin{Question}
\label{questaoinfbrp}
Is there an endomorphism  of a hyperbolic group with infinite image and infinite kernel for which the (synchronous) BRP holds? 
\end{Question} 

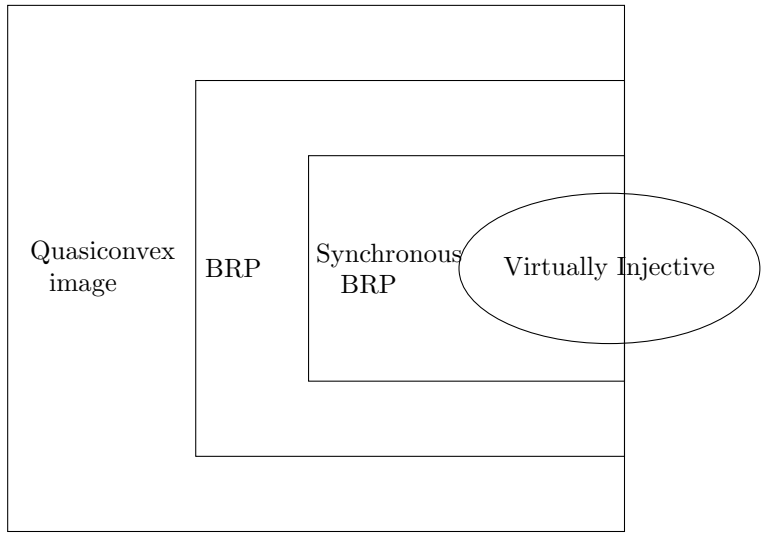
\begin{figure}[H]
\begin{center}
\begin{tikzpicture}
	\begin{scope} [fill opacity = 0.4]

    \draw[ draw = black] (-6.5,-2.5) rectangle (1.7,4.5);
        \draw[ draw = black] (-4,-1.5) rectangle (1.7,3.5);
\draw[  draw = black] (-2.5,-0.5) rectangle (1.7,2.5);

    \draw[ draw = black]  (1.5,1) ellipse (2 and 1);

 \node[opacity=1] at (-5.5,1) {\small\text{\parbox{1.4cm} {\centering Quasiconvex image}}};
       \node[opacity=1] at (-1.7,1)  {\small\text{\parbox{1.4cm} {\centering Synchronous BRP}}};

    \node[opacity=1] at (-3.5,1) {\small\text{BRP}};
        \node[opacity=1] at (1.5,1) {\small\text{Virtually Injective}};

    \end{scope}

\end{tikzpicture}
\caption{Nontrivial endomorphisms of hyperbolic groups}
\end{center}
\end{figure}

\subsection{Automatic groups}

In the case of automatic groups, we have essentially the same questions as in the case of hyperbolic groups except Question \ref{questaoinfbrp}, which we can answer affirmatively: consider $\Z\times \Z$, put $a=(1,0)$ and $b=(0,1)$ and take the structure $L=a^*b^*$. Let $\varphi\in \End(\Z\times \Z)$ defined by $(n,m)\mapsto (n,0)$. Then it is easy to see that the synchronous BRP holds for $(\varphi,L,L)$, but $\Ker(\varphi)\simeq \Z$ and Im($\varphi$)$\simeq \Z$.

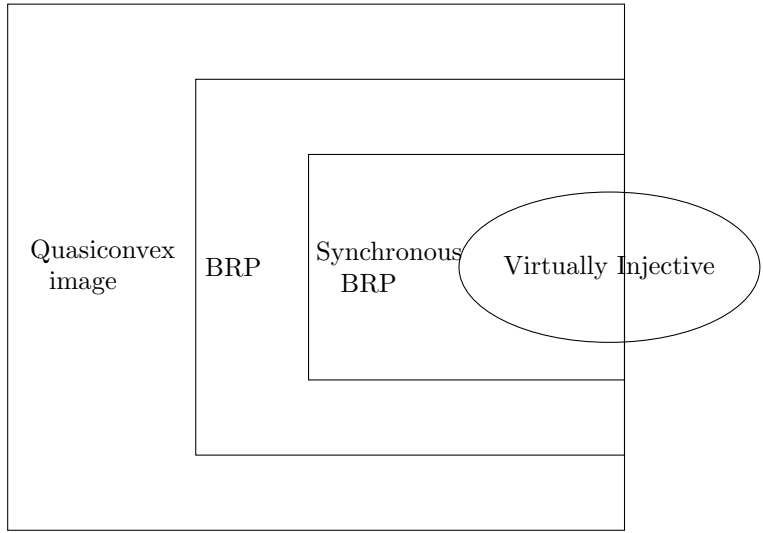
\begin{figure}[H]
\begin{center}
\begin{tikzpicture}
	\begin{scope} [fill opacity = 0.4]

    \draw[ draw = black] (-6.5,-2.5) rectangle (1.7,4.5);
        \draw[ draw = black] (-4,-1.5) rectangle (1.7,3.5);
\draw[  draw = black] (-2.5,-0.5) rectangle (1.7,2.5);

    \draw[ draw = black]  (1.5,1) ellipse (2 and 1);

 \node[opacity=1] at (-5.5,1) {\small\text{\parbox{1.4cm} {\centering Quasiconvex image}}};
       \node[opacity=1] at (-1.7,1)  {\small\text{\parbox{1.4cm} {\centering Synchronous BRP}}};

    \node[opacity=1] at (-3.5,1) {\small\text{BRP}};
        \node[opacity=1] at (1.5,1) {\small\text{Virtually Injective}};

    \end{scope}

\end{tikzpicture}
\caption{Nontrivial endomorphisms of automatic groups}
\end{center}
\end{figure}

\section{Further questions}
\label{secquestions}

The main questions arising from this work concern the satisfiability of the hypothesis in Corollary \ref{fix}.

\begin{Question}
In \cite{[GS91]}, the authors ask if a finite extension of a biautomatic group is also biautomatic. They also remark that to prove it, it suffices to show that the fixed subgroup of an automorphism of finite order is biautomatic. Can we use Theorem \ref{fix} to prove that, i.e., does such an endomorphism admit a language $L$ such that the synchronous BRP holds for $(\varphi, L, L)$? 
\end{Question}

\begin{Question}
Getting knowledge on the synchronous equivalence classes for some classes of automatic groups might be useful to prove fixed points results using Corollary \ref{fix}. There are infinitely many, but can we describe them in some sense?
\end{Question}

\section*{Acknowledgements}

The author is grateful to Pedro Silva for fruitful discussions of these topics, which greatly improved the paper.

The author was supported by the grant SFRH/BD/145313/2019 funded by Funda\c c\~ao para a Ci\^encia e a Tecnologia (FCT).

\end{document}